\documentclass{article}
\usepackage{authblk}
\usepackage{amsfonts, amssymb, amsmath, amsthm}
\usepackage{mathtools, bm, bbm, parskip}
\usepackage{graphicx,float,caption}
\usepackage{tabu, booktabs, tabularx}
\usepackage[shortlabels]{enumitem}
\usepackage{color}
\usepackage{xparse}
\usepackage{hyperref}
\usepackage[commandRef=Cref,createShortEnv, conf={one big link}]{proof-at-the-end}
\usepackage[noabbrev,capitalize]{cleveref}
\usepackage[maxnames=999,giveninits=true,style=numeric-comp]{biblatex}
\newcolumntype{Y}{>{\small\raggedleft\arraybackslash}X}

\newtheorem{theorem}{Theorem}
\newtheorem{definition}{Definition}

\newcommand{\indep}{\perp \!\!\! \perp}
\newcommand{\I}[1]{\mathbbm{1}{\left\{#1\right\}}}
\RenewDocumentCommand\P{mg}{\mathbb P\left(#1\IfNoValueF{#2}{\;\middle|\;#2}\right)}
\NewDocumentCommand\E{mg}{\mathbb E\left(#1\IfNoValueF{#2}{\;\middle|\;#2}\right)}
\NewDocumentCommand\Etm{mg}{\mathbb E_t\left(#1\IfNoValueF{#2}{\;\middle|\;#2}\right)}
\NewDocumentCommand\V{mg}{\mathrm{Var}\left(#1\IfNoValueF{#2}{\;\middle|\;#2}\right)}
\NewDocumentCommand\Vtm{mg}{\mathrm{Var}_t\left(#1\IfNoValueF{#2}{\;\middle|\;#2}\right)}
\newcommand{\HPi}{\left(\mathcal H_{\Pi}\right)}
\newcommand{\HC}{\left(\mathcal H_{\mathcal C}\right)}
\newcommand{\F}[1]{\mathcal F_{#1}}
\newcommand{\Ft}{\mathcal F_t}
\newcommand{\N}[1]{\mathcal{N}\left(#1\right)}
\newcommand{\tld}[1]{\widetilde{#1}}
\renewcommand{\hat}[1]{\widehat{#1}}
\makeatletter
\newcommand*\bigcdot{\mathpalette\bigcdot@{.5}}
\newcommand*\bigcdot@[2]{\mathbin{\vcenter{\hbox{\scalebox{#2}{$\m@th#1\bullet$}}}}}
\makeatother
\title{Non-parametric estimation of net survival under dependence between death causes}
\author[1]{Oskar Laverny}
\author[1]{Nathalie Grafféo}
\author[2]{Roch Giorgi}
\affil[1]{Aix Marseille Univ, INSERM, IRD, SESSTIM, Sciences Economiques \& Sociales de la Santé \& Traitement de l’Information Médicale, ISSPAM, Marseille, France.}
\affil[2]{Aix Marseille Univ, APHM, INSERM, IRD, SESSTIM, Sciences Economiques \& Sociales de la Santé \& Traitement de l’Information Médicale, ISSPAM, Hop Timone, BioSTIC, Biostatistique et Technologies de l’Information et de la Communication, Marseille, France.}
\begin{document}
\maketitle
\begin{abstract}
Relative survival methodology deals with a competing risks survival model where the cause of death is unknown.
This lack of information occurs regularly in population-based cancer studies.
Non-parametric estimation of the net survival is possible through the Pohar Perme estimator.
Derived similarly to Kaplan-Meier, it nevertheless relies on an untestable independence assumption.
We propose here to relax this assumption and provide a generalized non-parametric estimator that works for other dependence structures, by leveraging the underlying stochastic processes and martingales.
We formally derive asymptotics of this estimator, providing variance estimation and log-rank-type tests.
Our approach provides a new perspective on the Pohar Perme estimator and the acceptability of the underlying independence assumption.
We highlight the impact of this dependence structure assumption on simulation studies, and illustrate them through an application on registry data relative to colorectal cancer, before discussing potential extensions of our methodology.
\end{abstract}

\section{Introduction}

Survival analysis deals with the modeling of time-to-event data, in particular under censorship~\cite{andersenStatisticalModelsBased1993,fleming2013counting,Collett2023}.
The tools developed in this field are nowadays crucial in medical research, when investigating the impact of certain diagnosis and/or treatment of diseases on the patient's overall survival.
The main practical output is the statistical estimation of the death probability across time, possibly conditionally on covariates~\cite{Jenkins2005,Survival2010}.
The theory also provides analytical tools (test statistics, asymptotics, etc.) to help practitioners analyze their datasets and identify important variables affecting the survival of their patients.
Not technically restricted to human lifetime, survival analysis is a fundamental tool in many fields, from operational research and actuarial sciences to chemistry, biostatistics, and of course medical science.

In population-based cancer studies, the reported cause of death -- assumed binary, death from the studied cancer or from other causes -- is usually unavailable or unreliable in the available datasets.
The relative survival methodology \cite{poharpermeEstimationRelativeSurvival2012, jooste2013unbiased} takes this particularity into account to evaluate the excess mortality -- due to cancer -- by assuming the general population mortality rates to be known, extracted from population life tables.

The methodology allows epidemiologists to estimate an excess mortality with respect to cancer only.
It also allows comparing the excess mortality of groups with different populational distributions (e.g., for spatial, temporal, or age disparities within and between cohorts) with proper statistical tools (asymptotic results, tests, etc.).
Many national and international health agencies use it to report cancer statistics: the World Health Organization (WHO), the International Agency for Research on Cancer (IARC), the FRANCIM network in France~\cite{bossard2007survival,cowppli2017survival}, or international research programs EUROCARE~\cite{berrino2007survival, brenner2012progress} at European level or CONCORD~\cite{coleman2008cancer, allemani2018global}.
They provide excess mortality estimates in their periodic reports and publications, do cross-country comparisons, trend analysis, and highlight potential areas of improvement for public health decision-making.
By overcoming unreliable cause-of-death reportings, the methodology offers practical insights for equity, policy, and clinical cares. 

Let us now describe a relative survival problem formally, following e.g., \cite{esteveRelativeSurvivalEstimation1990}.
Consider for this three random times $E,P$ and $C$, respectively the times, from diagnosis, to death from the studied cause (E stands for "Excess"), to death from other causes (P stands for "Population") and to censorship.
Denote the time to death from any cause by $O = E \wedge P$ (O stands for "Overall", while the symbol $\wedge$ denotes the minimum).
In relative survival setting, times $E,P,O$ and the cause of death indicatrix $\I{E \le P}$ are not observed: only the time to event $T = O \wedge C$ and the censoring indicatrix $\Delta = \I{T \le C}$ are observed, together with a vector of covariates $\bm X$.
This separates relative survival settings from standard competing risks setting (where $\I{E \le P}$ would also be observed).
Note that we cannot distinguish a death due to the studied cancer from a death due to other causes.

A first standard assumption in the net survival literature is that the distribution of $P \;|\; \bm X$ (other causes of death, differentiated along the covariates) is known and can be derived from a reference population life table for each individual.
This requires assuming that the prevalence of deaths due to the interesting cause (the excess mortality) in these life tables is sufficiently small to be negligeable.
To simplify the exposition, assume that the rest of the model (that is the distributions of $E$, $C$, and the dependence structure between $(E,P,C)$) does not depend on covariates $\bm X$.
Relative survival's goal is to estimate the distribution of $E$, and this distribution only: distributions of $P$, $C$, and the dependence structure of $(E,P,C)$ are considered nuisances.

Let $(\bm X_i, P_i, E_i, C_i, O_i = P_i \wedge E_i, T_i = O_i \wedge C_i, \Delta_i = \I{T_i \le C_i})_{i \in 1,...,n}$ be an $i.i.d$ $n$-sample of these random variables.
Note that only $(\bm X_i, T_i, \Delta_i)$ are observed: we thus let $\left(\Omega,\mathcal A,\left\{\Ft, t\in \mathbb R_+\right\},\mathbb P \right)$ be the associated filtered probability space, where the filtration corresponds to our observations: 
$$\Ft = \sigma\left\{\bm{X}_i,\left(T_i,\Delta_i\right): T_i \le t,\;\forall i \in 1,..,n\right\}.$$

Survival analysis usually considers the different times, in particular $T$ and $C$, to be independent.
Relaxation of this hypothesis is an active field of research, see \cite{czadoDependentCensoringBased2021,deresaCopulaBasedCox2023,deresaCopulabasedInferenceBivariate2022,deresaMultivariateNormalRegression2020,deresaSemiparametricModellingEstimation2021} for recent developments.
Relative survival analysis is no outlier and commonly assumes that $E,P$ and $C$ are mutually independent.
If the epidemiological interpretation makes this assumption hard to justify -- the reasons causing death by the studied disease are related to the reasons causing other deaths, and even to the reasons causing loss of follow-up and therefore censorship -- it is still very practical and simplifies drastically the estimation process.
This assumption is therefore central in the literature, and this is the assumption we want to challenge here in the case of relative survival.

We consider in this paper a relaxation of the independence assumption that only assumes independent censorship, of the form $(E,P) \indep C$, while the dependence structure of $(E,P)$ has to be taken into account.
We express this dependence structure through the survival copula $\mathcal{C}$ of the random vector $(E,P)$ (see \cite{nelsenIntroductionCopulas2006} or more recently \cite{durante2015a}).
Our assumption writes:
\begin{equation}\label{eq:first_equation}
\HC:\;S_{E \wedge P}(t) = \mathcal{C}\left(S_{E}(t),S_{P}(t)\right).
\end{equation} 

In particular, denoting $\Pi(u_1,u_2) = u_1u_2$ the bivariate independence copula, the standard assumption writes simply $\HPi$.
This dissection of the dependence structure of $(E,P,C)$ is crucial to match closely the nature of the observations: we have more statistical information on the joint distribution of $(O,C)$ (recall $O = E\wedge P$) than on the joint distribution of $(E,P)$.
More precisely, the joint distribution of $(O, C)$ is tied to the observations $\left(T_i,\Delta_i\right)_{i \in 1,...,n}$ exactly as in a classical censored survival problem.
The supposed independence of the censorship is there only to simplify our expressions.
On the other hand, the dependence structure of $(E,P)$ is completely unattainable from our observations, and we must rely on an assumption $\HC$ -- but any copula will do.

There are a few standard non-parametric estimators of the distribution of $E$ under $\HPi$ in the literature, starting in the 1960s with the seminal work of Ederer~\cite{ederer1961relative,ederer1959effect}, Hakulinen~\cite{hakulinen1982cancer} in the 1980s, and more recently Pohar Perme~\cite{poharpermeEstimationRelativeSurvival2012} in 2012.
However, groundings and justifications of these estimators is a bit blurry.
Under $\HC$ on the other hand, the literature does not provide non-parametric estimators, and the parametric estimation routines from \cite{adatorwovorParametricApproachRelaxing2021} are quite recent.

We propose here a generalized version of the non-parametric Pohar Perme estimator from $\HPi$ to $\HC$.
Leveraging tools from counting processes theory, we formally derive asymptotics of this new estimator in \cref{sec:principal}, providing variance estimation and associated log-rank-type tests under $\HC$.
This methodological breakthrough provides a new perspective on the Pohar Perme estimator and the acceptability of the underlying independence assumption.
\Cref{sec:simus} contains simulations studies to understand the tradeoff of our approach and evaluate performance of the proposed schemes, while \cref{sec:examples} gives a real data example on colorectal cancer.
Both sections highlight the importance of the choice of the copula on the results, and \cref{sec:conclusion} concludes.

\section[Relative survival under (HC)]{Relative survival under $\HC$}\label{sec:principal}

In this section, we formally derive non-parametric estimators for the excess hazard, its variance, associated asymptotics and log-rank-type tests.
Proofs and secondary technical statements are differed to \cref{apx:proofs} to ensure a short and clear exposition of the main results.

\subsection{Estimation of the excess hazard}

We derive first a non-parametric estimator of the net survival under $\HC$.
The proposed derivation, when applied to the particular case of $\HPi$, coherently ends up on the Pohar Perme estimator~\cite{poharpermeEstimationRelativeSurvival2012}, the \emph{de facto} standard non-parametric estimator of the net survival under $\HPi$.
Our derivation thus provides a new perspective on the Pohar Perme estimator and its weighting scheme, clarifying how it relies on the hypothesis $\HPi$.

Define first the following stochastic processes: 
\begin{alignat*}{3}
    N(t) &= \I{O \le t, O \le C} &&\;\textit{ (Uncensored deaths process)}\\
    Y(t) &= \I{O \ge t, C \ge t} &&\;\textit{ (At-risk process)}\\
    N_E(t) &= \I{E \le t, E \le C}&&\;\textit{ (Excess uncensored deaths process)}\\
    Y_E(t) &= \I{E \ge t, C \ge t}&&\;\textit{ (Excess at-risk process)}\\
\end{alignat*}
We similarly define individual processes $N_i,Y_i, N_{E,i}$ and $Y_{E,i}$ for all (independent) patients $i \in 1,..,n$.
Note that our filtration can be re-written in terms of $N_i,Y_i$ as follows:
$$\Ft = \sigma\left\{ \bm X_i, (N_i(s),Y_i(s)): 0 \le s \le t, i \in 1,...,n\right\}.$$
We also denote $\Etm{\cdot}$ as a shorthand for $\E{\cdot}{\F{t^-}}$.
From the standard Doob-Meyer decomposition of the process $N_i(t)$, we know that, first, the process $M_i(t)$ is a local square integrable $\Ft$-martingale, and, second, that $Y_i(t)$ is a bounded $\Ft$-predictable process.

The central issue in relative survival is the construction of estimators of the excess mortality distribution, through the excess survival function $S_E(t)$ and/or the excess hazard function $\Lambda_E(t) = -\ln S_E(t)$.
Indeed, such estimators are usually based on $\left(N_{E,i},Y_{E,i}\right)_{i \in 1,...,n}$, whereas here only $\left(N_{i},Y_{i}\right)_{i \in 1,...,n}$ are observed.
We thankfully found a formal link between the two.
To express it, define first $a,b,c$ and $a_i,b_i,c_i$ as:
\begin{align*}
    a(t) &= \P{P \ge t}{E = t}   \;\text{ and }\; a_i(t) = \P{P_i \ge t}{E_i = t},\\
    b(t) &= \P{P = t}{E \ge t}   \;\text{ and }\; b_i(t) = \P{P_i = t}{E_i \ge t},\\
    c(t) &= \P{P \ge t}{E \ge t} \;\text{ and }\; c_i(t) = \P{P_i \ge t}{E_i \ge t}.
\end{align*}

Denote finally $M(t) = N(t) - \int_0^t Y(s)\partial\Lambda_O(s)$ the martingale classically associated with the process $N$, and $M_i$ similarly.

\begin{lemmaE}[Excess survival processes][end, restate]\label{lem:tech}
    The following points hold: 
    \begin{enumerate}[(i)]
        \item\label{lem:tech:1} $\partial N_E(t) = \frac{1}{a(t)} \E{\partial N(t)}{E,C}- \frac{b(t)}{a(t)c(t)} \E{Y(t)}{E, C}$
        \item\label{lem:tech:Y} $Y_E(t) = \frac{1}{c(t)} \E{Y(t)}{E, C}$
        \item\label{lem:tech:2} $\partial \Lambda_E(t) = \frac{c(t)}{a(t)}\left(\partial \Lambda_O(t) - \frac{b(t)}{c(t)}\right).$
        \item\label{lem:tech:3} The process $N_E$ admits the following Doob-Meyer decomposition: 
        $$\partial N_E(t) = \partial M_E(t) + Y_E(t) \partial \Lambda_E(t),$$
        where \begin{itemize}
            \item $Y_E(t) \partial \Lambda_E(t)$ is $\Ft$-predictable,
            \item $M_E(t) = \E{\int_0^t \frac{\partial M(s)}{a(s)}}{E,C}$ is a $\Ft$-martingale.
        \end{itemize}
    \end{enumerate}
\end{lemmaE}
\begin{proofE}
    We prove items in the exposing order.
    Remark first that$\P{P=E=t}=0$ by absolute continuity, and then: \begin{align*}
        \E{\partial N(t)}{E,C} 
            &= \P{E \wedge P = t, t \le C}{E,C}\\
            &= \P{P \ge t, E = t, t \le C}{E,C} + \P{t \le E, P = t, t \le C}{E,C}\\
            &= \P{P \ge t}{E = t}\I{E = t, E \le C} + \P{P = t}{E \ge t}\I{E \ge t, C \ge t}\\
            &= a(t)  \partial N_E(t) + b(t) Y_E(t).
    \end{align*}
    Moreover,
    \begin{align*}
        \E{Y(t)}{E,C}
            &= \P{E\wedge P \ge t, C \ge t}{E,C}\\
            &= \P{P \ge E \ge t, C \ge t}{E,C} + \P{E \ge P \ge t, C \ge t}{E,C}\\
            &= \P{P \ge E}{E \ge t} Y_E(t) + \P{E \ge P \ge t}{E \ge t}Y_E(t)\\
            &= \P{P \ge t}{E \ge t} Y_E(t)\\
            &= c(t) Y_E(t).
    \end{align*}
    \Cref{lem:tech:1,lem:tech:Y} follow from inverting this system of two equations to extract $(\partial N_E(t), Y_E(t))$.
    To prove \cref{lem:tech:2}, recall that $\partial \Lambda_X(t) = \frac{\P{X = t}}{\P{X \ge t}}$ for $X \in {E,O}$, and expand $a(t),b(t),c(t)$:
    \begin{align*}
        \frac{c(t)}{a(t)}\left(\partial \Lambda_O(t) - \frac{b(t)}{c(t)}\right)
        &= \frac{\P{P \ge t}{E \ge t}}{\P{P \ge t}{E = t}}\left(\frac{\P{O = t}}{\P{O \ge t}} - \frac{\P{P = t}{E \ge t}}{\P{P \ge t}{E \ge t}}\right)\\
        &= \frac{\P{E = t}}{\P{E \ge t}}\frac{\P{P \ge t , E \ge t}}{\P{P \ge t , E = t}}\left(\frac{\P{O = t}}{\P{O \ge t}} - \frac{\P{P = t , E \ge t}}{\P{P \ge t , E \ge t}}\right)\\
        &= \partial \Lambda_E(t)\frac{\P{O \ge t}}{\P{P \ge t , E = t}}\left(\frac{\P{O = t} - \P{P = t , E \ge t}}{\P{O \ge t}}\right)\\
        &= \partial \Lambda_E(t)\frac{\P{O = t} - \P{P = t , E \ge t}}{\P{P \ge t , E = t}}\\
        &= \partial \Lambda_E(t).
    \end{align*}

    To prove \cref{lem:tech:3}, using \cref{lem:tech:1,lem:tech:Y,lem:tech:2}, consider first that 
    \begin{align*}
        \Etm{\partial N_E(t)}
        &= \Etm{\frac{1}{a(t)}\E{\partial N(t)}{E,C} - \frac{b(t)}{a(t)c(t)} \E{Y(t)}{E,C}}\\
        &= \E{\frac{1}{a(t)}\Etm{\partial N(t)} - \frac{b(t)}{a(t)c(t)} Y(t)}{E,C}\\
        &= \E{\frac{Y(t)}{a(t)}\left(\partial\Lambda_O(t) - \frac{b(t)}{c(t)}\right)}{E,C}\\
        &= \E{Y_E(t) \partial \Lambda_E(t)}{E,C}\\
        &= Y_E(t) \partial \Lambda_E(t).
    \end{align*}
    Thus,
    \begin{align*}
        \partial M_E(t) &= \partial N_E(t) - Y_E(t) \partial \Lambda_E(t)\\
        &= \E{\frac{\partial N(t)}{a(t)} - \frac{b(t) Y(t)}{a(t)c(t)} - \frac{Y(t)}{c(t)}\frac{c(t)}{a(t)}\left(\partial\Lambda_O(t) - \frac{b(t)}{c(t)}\right)}{E,C}\\
        &= \E{\frac{\partial N(t)}{a(t)} - \frac{Y(t)\partial\Lambda_O(t)}{a(t)}}{E,C}\\
        &= \E{\frac{\partial M(t)}{a(t)}}{E,C}.
    \end{align*}
    So that $M_E(t) = \E{\int_0^t \frac{\partial M(s)}{a(s)}}{E,C}$ as wanted.
    Since the function $a$ is predictable, $M_E$ is a martingale and the decomposition is unique.
\end{proofE}
Since \cref{lem:tech} hold for each one of our independent and identically distributed patients, a natural estimator of $\partial \Lambda_E$ could be constructed alike Kaplan-Meier~\cite{kaplanNonparametricEstimationIncomplete1958} as:
\begin{equation*}
    \frac{\frac{1}{n} \sum_{i=1}^n \partial N_{E,i}(t)}{\frac{1}{n} \sum_{i=1}^n Y_{E,i}(t)}.
\end{equation*}
However, $\partial N_{E,i}(t)$ and $Y_{E,i}(t)$ are not observable and need to be estimated.
To obtain observables, we propose to simply drop the conditional expectations operators from \cref{lem:tech:1,lem:tech:Y} in \cref{lem:tech}, that is to replace conditional expectations by their stochastic unconditional counterpart.
This process leads to the definition of a relevant estimator of $\partial\Lambda_E$ as: 
\begin{align}
    \partial \tld{N}_{E,i}(t) &= \frac{\partial N_i(t)}{a_i(t)} - \frac{b_i(t)}{a_i(t)c_i(t)}Y_i(t)\label{eq:def_dntld}\\
    \tld{Y}_{E,i}(t) &= \frac{Y_i(t)}{c_i(t)}\nonumber\\
    \partial \tld{\Lambda}_E (t) &= \frac{\sum_{i=1}^n \partial \tld{N}_{E,i}(t)}{\sum_{i=1}^n \tld{Y}_{E,i}(t)}\nonumber.
\end{align}

Focus now on the coefficients $a,b,c$.
They can be expressed in terms of the copula and survival functions of $E$ and $P$. 
Indeed, assuming $(E,P)$ is an absolutely continuous random vector, let us denote the partial derivatives of the copula $\mathcal C$ as $\mathcal C_j(\bm u) = \frac{\partial \mathcal C}{\partial u_j}(\bm u), j \in {1,2}$. 
We then have:
\begin{align*}
    a(t) &= \mathcal C_1\left(S_E(t), S_P(t)\right) \\
    b(t) &= -\mathcal C_2\left(S_E(t), S_P(t)\right) \frac{\partial S_P(t)}{S_{E}(t)}\\
    c(t) &= \frac{\mathcal C(S_{E}(t), S_{P}(t))}{S_{E}(t)} = \frac{S_O(t)}{S_E(t)}.
\end{align*}
Under $\HPi$, we remark that $a(t) = c(t) = S_P(t)$ and $b(t) = -\partial S_P(t)$, which do not depend on the excess distribution, and are hence observables.
$\partial \tld{\Lambda}_E$ is thus observable, and is in fact already known as being the Pohar Perme estimator~\cite{poharpermeEstimationRelativeSurvival2012}.

However, under $\HC$, $a(t),b(t)$ and $c(t)$ depend on the (unknown) $S_E(t)$.
The individual coefficients $a_i,b_i$ and $c_i$ thus need to be estimated alongside $\partial\Lambda_E$, which leads to \cref{def:gen_ppe}.

\begin{definition}[Generalized Pohar Perme estimator]\label{def:gen_ppe}
    We call \emph{generalized Pohar Perme estimator} the solution $\hat{\Lambda}_E$ of the differential equation
    \begin{equation}\label{eq:gen_ppe}
        \partial\hat{\Lambda}_E(t) =  \frac{\sum_{i=1}^n  \partial \hat{N}_{E,i}(t)}{\sum_{i=1}^n \hat{Y}_{E,i}(t)},
    \end{equation}
    where $\hat{S}_E(t) = \exp\left\{-\hat{\Lambda}_E(t)\right\}$ and for all $i \in 1,...,n$,
    \begin{align*}
        \partial \hat{N}_{E,i}(t) &= \frac{\partial N_i(t)}{\hat{a}_i(t)}  - \frac{\hat{b}_i(t)Y_i(t)}{\hat{a}_i(t)\hat{c}_i(t)},\\
        \hat{Y}_{E,i}(t) &= \frac{Y_i(t)}{\hat{c}_i(t)},\\
        \hat{a}_i(t) &= \mathcal C_1\left(\hat{S}_E(t),S_{P_i}(t)\right),\\
        \hat{b}_i(t) &= \mathcal C_2\left(\hat{S}_E(t), S_{P_i}(t)\right) \frac{-\partial S_{P_i}(t)}{\hat{S}_E(t)},\\
        \hat{c}_i(t) &= \frac{\mathcal C\left(\hat{S}_E(t), S_{P_i}(t)\right)}{\hat{S}_E(t)}.
    \end{align*}
\end{definition}

Under $\HPi$, we have $\hat{\Lambda}_E = \tld{\Lambda}_E$.
In all generality however, the plug-in of $\hat{\Lambda}_E(t)$ in the right-hand-side of the (unseparable) differential equation is necessary, and a non-linear equation in $\partial \hat{\Lambda}_E(t)$ needs to be solved at each time step.
As previous estimators under $\HPi$, the obtained $\partial\hat{\Lambda}_E$ process is piecewise continuous, with jumps at event times $T_1,...,T_n$, and always negative (except at jump points).
Solving the differential equation must therefore use a very dense mesh $t_1,...,t_N$ that includes observed failure times $T_1,...T_n$.
These characteristics were already present under $\HPi$, and thus our implementation mirrors the discretization scheme from~\cite{poharpermeEstimationRelativeSurvival2012} (implicit Euler).

Remark that our construction would allow for different specifications of the dependence structure for each individual, but studying this case is beyond the scope of this paper.

\subsection{Variance estimation}

We now focus on the estimation of the variances of $\tld{\Lambda}_E$ and $\hat{\Lambda}_E$.
Let us denote $$\Xi(t) = \int_{0}^t \frac{\sum_{i=1}^n \frac{\partial M_i(s)}{a_i(s)}}{\sum_{i=1}^n \frac{Y_i(s)}{c_i(s)}}.$$

\begin{lemmaE}[$\tld{\Lambda}_E$'s Doob-Meyer decomposition][end,restate]\label{lem:DMtldlambda}
    $\Xi(t)$ is a local square integrable martingale, and the Doob-Meyer decompositions of $\tld{\Lambda}_E$ writes $$\tld{\Lambda}_E(t) = \Lambda_E(t) + \Xi(t).$$
\end{lemmaE}
\begin{proofE}
    Let us denote $\forall i \in 1,...,n$:
    \begin{align*}
        \partial\tld{M}_{E}(t) &= \frac{\partial M(t)}{a(t)}\\
        \partial\tld{M}_{E,i} &= \frac{\partial M_i(t)}{a_i(t)}\\
        \partial \Xi(t) &= \frac{\sum_{i=1}^n \partial\tld{M}_{E,i}(t)}{\sum_{i=1}^n \tld{Y}_{E,i}(t)}.
    \end{align*}
    We will in fact prove that both $\tld{M}_{E,i}$ and $\Xi$ are local square integrable martingales, and that we have the two following Doob-Meyer decompositions: 
    \begin{itemize}
        \item $\forall i,$ $\partial \tld{N}_{E,i}(t) = \partial\tld{M}_{E,i}(t) + \tld{Y}_{E,i}(t) \partial \Lambda_E(t),$
        \item $\tld{\Lambda}_E(t) = \Lambda_E(t) + \Xi(t).$
    \end{itemize}

    Recall first that, for all $i \in 1,...,n$, $$\partial \tld{N}_{E,i}(t) =  \frac{1}{a_i(t)} \partial N_i(t) - \frac{b_i(t)}{a_i(t)c_i(t)} Y_i(t).$$
    Since $\Etm{\partial N_i(t)} = Y_i(t)\partial \Lambda_{O_i}(t)$, we have by \cref{lem:tech:2} in \cref{lem:tech} that:
    \begin{align*}
        \Etm{\partial \tld{N}_{E,i}(t)} 
        &= \frac{1}{a_i(t)} Y_i(t)\partial \Lambda_{O_i}(t) - \frac{b_i(t)}{a_i(t)c_i(t)} Y_i(t)\\
        &= \frac{Y_i(t)}{a_i(t)}\left(\partial \Lambda_{O_i}(t) - \frac{b_i(t)}{c_i(t)}\right)\\
        &= \frac{Y_i(t)}{c_i(t)}\partial \Lambda_{E_i}(t)\\
        &= \tld{Y}_{E,i}(t) \partial \Lambda_E(t).
    \end{align*}
    Subtracting the two gives $\partial \tld{M}_{E,i}(t) = \frac{\partial M_i(t)}{a_i(t)}$ as required, and since $
    \tld{M}_{E,i}(t)$ is a local square integrable martingale, this decomposition is unique.

    Then, since $a_i(t),b_i(t),c_i(t)$ are constants and since the processes $Y_i$ are $\Ft$-predictables, we have that
    \begin{align*}
        \Etm{\partial \tld{\Lambda}_E(t)} 
        &= \Etm{\frac{\sum_{i=1}^n \partial \tld{N}_{E,i}(t)}{\sum_{i=1}^n \tld{Y}_{E,i}(t)}}\\
        &=\frac{\sum_{i=1}^n \Etm{\partial \tld{N}_{E,i}(t)}}{\sum_{i=1}^n \tld{Y}_{E,i}(t)}\\
        &=\frac{\sum_{i=1}^n \tld{Y}_{E,i}(t) \partial \Lambda_E(t)}{\sum_{i=1}^n \tld{Y}_{E,i}(t)}\\
        &= \partial \Lambda_E(t)
    \end{align*}
    Subtracting, we obtain the required expression for $\partial\Xi(t)$ as follows:
    $$\partial\Xi(t) = \partial \tld{\Lambda}_E(t) - \partial \Lambda_E(t) = \frac{\sum_{i=1}^n \frac{1}{a_i(t)}\partial M_i(t)}{\sum_{i=1}^n \frac{Y_i(t)}{c_i(t)}}.$$
    Since $Y_i$ is predictable, again, $\Xi$ is a local square integrable martingale and this decomposition is unique.
\end{proofE}
For a process $X(t)$, we denote by $\left[X\right](t)$ its associated optional process (see \cite{andersenStatisticalModelsBased1993}).
A straightforward consequence of \cref{lem:DMtldlambda} is then given by \cref{prop:variance}.

\begin{propositionE}[Variance of $\tld{\Lambda}_E(t)$][end, restate]\label{prop:variance} We have
    $$\V{\tld{\Lambda}_E(t)} = \E{\left[\Xi\right](t)} \text{ and } \left[\Xi\right](t) = \int_{0}^t \frac{\sum_{i=1}^n \frac{\partial N_i(s)}{a_i(s)^2}}{\left(\sum_{i=1}^n \frac{Y_i(s)}{c_i(s)}\right)^2}.$$
\end{propositionE}
\begin{proofE}
    From \cite[formula (2.47) and (2.49) p. 56]{aalenSurvivalEventHistory2008}, denoting $\left[\Xi\right]$ denote the optional processes of $\Xi$, \cref{lem:DMtldlambda} implies that $$\V{\tld{\Lambda}_E(t)} = \E{\left[\Xi\right](t)}.$$
    Due to orthogonality of martingales $M_1,...,M_n$, and since $\forall i$, $$\frac{\frac{1}{a_i(s)}}{\sum_{i=1}^n \frac{Y_i(s)}{c_i(s)}}$$ are predictable processes, we complete the argument by:
    \begin{align*}
        \partial \left[\Xi\right](t) &= \partial \left[\int_{0}^t \frac{\sum_{i=1}^n \frac{1}{a_i(s)}\partial M_i(s)}{\sum_{i=1}^n \frac{Y_i(s)}{c_i(s)}}\right](t)\\
        &= \partial \left[\sum_{i=1}^n \int_{0}^t \frac{\frac{1}{a_i(s)}}{\sum_{i=1}^n \frac{Y_i(s)}{c_i(s)}} \partial M_i(s)\right](t)\\
        &= \partial \sum_{i=1}^n \int_{0}^t \left(\frac{\frac{1}{a_i(s)}}{\sum_{i=1}^n \frac{Y_i(s)}{c_i(s)}}\right)^2 \partial N_i(s)\\
        &= \frac{\sum_{i=1}^n \frac{1}{a_i(t)^2}\partial N_i(t)}{\left(\sum_{i=1}^n \frac{Y_i(t)}{c_i(t)}\right)^2}.
    \end{align*}
\end{proofE}
A good estimator for the variance of $\tld{\Lambda}_E(t)$ is thus simply $$\tld{\sigma}_E^2(t) =  \left[\Xi\right](t).$$
Note that this estimator is unbiased by construction.
Under $\HPi$, $\tld{\sigma}_E^2(t)$ is feasible (relies only on observable quantities), and corresponds to the variance already obtained in~\cite{poharpermeEstimationRelativeSurvival2012}.
However, under $\HC$, once again $\tld{\sigma}_E^2(t)$ is not feasible, and thus we propose to use the straightforward plug-in estimator $\hat{\sigma}_E^2(t)$:

\begin{equation}\label{eq:gen_var}
    \hat{\sigma}_E^2(t) = \int_{0}^t \frac{\sum_{i=1}^n \frac{\partial N_i(s)}{\hat{a}_i(s)^2}}{\left(\sum_{i=1}^n \frac{Y_i(s)}{\hat{c}_i(s)} \right)^2}.
\end{equation}

The performance of this last plug-in will be discussed through experiments in \cref{sec:simus,sec:examples}.

\subsection{Log-rank type test}

Suppose that the patients are dispersed in disjoint groups forming a partition $G = \left\{g_1,..,g_r\right\}$ of the set of indices $\left\{1,..,n\right\}$.
We want to construct a test for the equality of the excess mortality distributions in these groups.
That is, under the weaker assumption that $E_1,...,E_n$ are only identically distributed inside each of the groups, we want to test the assumption that $E_1,...E_n$ are, in fact, identically distributed among our population.
Since the cumulative hazard function characterizes the distribution of a positive random variable, this hypothesis can be formalized as follows: 
$$\left(H_0\right): \forall g \in G, \forall i \in g, \;\Lambda_{E_i} = \Lambda_{E}.$$

Let us denote $\tld{Y}_{E,g} = \sum_{i \in g} \tld{Y}_{E,i}$ for any group $g \in G$, and $\tld{Y}_{E,\bigcdot} = \sum_{g \in G} \tld{Y}_{E,g}$.
Similarly, denote  $\tld{N}_{E,g} = \sum_{i \in g} \tld{N}_{E,i}$ and $\tld{N}_{E,\bigcdot} = \sum_{g \in G} \tld{N}_{E,g}$.

Define finally the vectors $\bm R(t),\bm Z(t)$, the matrix $\bm \Gamma(t)$ and the test statistic $\tld{\chi}(T)$ by: 
\begin{align*}
    R_g(t) &= \frac{\tld{Y}_{E,g}(t)}{\tld{Y}_{E,\bigcdot}(t)}\\
    Z_{g}(t) &= \tld{N}_{E,g}(t) - \int_{0}^t R_g(s) \partial \tld{N}_{E,\bigcdot}(s)\\
    \Gamma_{g,h}(t) &= \sum_{\ell\in G} \int_{0}^t \left(\delta_{\ell,g} - R_g(s)\right)\left(\delta_{\ell,h} - R_h(s)\right) \sum_{i \in \ell} \frac{\partial N_i(s)}{a_i(s)^2}.\\
    \tld{\chi}(T) &= \bm Z(T)'\bm \Gamma(T)^{-1} \bm Z(T)
\end{align*}

\begin{lemmaE}[Properties of $\bm R, \bm Z$ and $\bm \Gamma$][all end]\label{lem:test}
    Let $T < \infty$.
    Under $\left(H_0\right)$, assuming that there exists an $\epsilon > 0:\;a_i(t) > \epsilon,c_i(t)>\epsilon$ over $t \in [0,T]$, the following points hold over $t \in [0,T]$,
    \begin{enumerate}[(i)]
        \item\label{itm:Z_centered_lsim} $\bm Z$ is a centered local square integrable martingale
        \item\label{itm:Cov_Z} $\mathrm{Cov}(\bm Z(t)) = \E{\bm \Gamma(t)}$
        \item\label{itm:cv_Gamma} $n^{-1} \bm \Gamma(t) \xrightarrow[n \to \infty]{\mathbb P} \bm V(t)$, $\bm V$ is deterministic, and both $\bm \Gamma(t)$ and $\bm V(t)$ are semi-definite positives.
        \item\label{itm:cv_Z} $n^{-\frac{1}{2}}\bm Z(t) \xrightarrow[n\to\infty]{\mathcal D} \N{0,\bm V(t)}$
        \item\label{itm:rank_V} $\mathrm{Ker}(\bm V(t)) = \mathrm{Vect}(\bm 1)$
    \end{enumerate}
\end{lemmaE}
\begin{proofE}
    Consider throughout the proof that $\left(H_0\right)$ holds, and recall definitions of vectors $\bm R(t),\bm Z(t)$, the matrix $\bm \Gamma(t)$ and the test statistic $\tld{\chi}(T)$:
    \begin{align*}
        R_g(t) &= \frac{\tld{Y}_{E,g}(t)}{\tld{Y}_{E,\bigcdot}(t)}\\
        Z_{g}(t) &= \tld{N}_{E,g}(t) - \int_{0}^t R_g(s) \partial \tld{N}_{E,\bigcdot}(s)\\
        \Gamma_{g,h}(t) &= \sum_{\ell\in G} \int_{0}^t \left(\delta_{\ell,g} - R_g(s)\right)\left(\delta_{\ell,h} - R_h(s)\right) \sum_{i \in \ell} \frac{\partial N_i(s)}{a_i(s)^2}.\\
        \tld{\chi}(T) &= \bm Z(T)'\bm \Gamma(T)^{-1} \bm Z(T)
    \end{align*}
    Let us prove points in their exposing order.
    
    \emph{Proof of \cref{itm:Z_centered_lsim}:} From \cref{lem:DMtldlambda}, we have, $\forall g \in G \cup \{\bigcdot\}$:
    $$\partial \tld{N}_{E,g}(t) = \partial \tld{M}_{E,g}(t) + \tld{Y}_{E,g}(t) \partial \Lambda_{E}(t).$$
    Thus:
    \begin{align}
        \partial Z_{g}(t) &= \partial \tld{N}_{E,g}(t) - R_g(t)\partial \tld{N}_{E,\bigcdot}(t)\nonumber\\
        &= \partial \tld{M}_{E,g}(t) + \tld{Y}_{E,g}(t) \partial \Lambda_E(t) - R_g(t) \left(\partial \tld{M}_{E,\bigcdot}(t) + \tld{Y}_{E,\bigcdot}(t) \partial \Lambda_E(t)\right)\nonumber\\
        &= \left(\partial \tld{M}_{E,g}(t) - R_g(t)\partial \tld{M}_{E,\bigcdot}(t)\right) + \left(\tld{Y}_{E,g}(t) \partial \Lambda_E(t) - R_g(t)\tld{Y}_{E,\bigcdot}(t) \partial \Lambda_E(t) \right)\nonumber\\
        &= \sum_{\ell \in G} \left(\delta_{g,\ell} - R_g(t)\right)\partial \tld{M}_{E,\ell}(t) + \tld{Y}_{E,g}(t) \left(\partial \Lambda_E(t) - \partial \Lambda_E(t)\right)\nonumber\\
        &= \sum_{\ell \in G} \left(\delta_{g,\ell} - R_g(t)\right) \partial \tld{M}_{E,\ell}(t)\label{eq:Z}
    \end{align}

    From now on, we denote 
    $$H_{g,\ell}(t) = \delta_{g,\ell} - R_g(t).$$

    Since $\left(\tld{M}_{E,g}\right)_{g\in G}$ are uncorrelated local square integrable martingales, and $\left(H_{g,\ell}\right)_{g,\ell\in G}$ are predictable processes, $Z_g$ clearly are centered local square integrable martingales themselves.

    \emph{Proof of \cref{itm:Cov_Z}:} We want to show that $\mathrm{Cov}(\bm Z(t)) = \E{\bm \Gamma(t)}$.
    For that, it would suffice to show that $\forall g,h\in G,\; \Gamma_{g,h}(t) = \left[Z_g,Z_h\right](t)$, so we now compute $\left[Z_g,Z_h\right](t)$.

    Since $\left(H_{g,\ell}\right)_{g,\ell\in G}$ are predictable processes and $\left(\tld{M}_{E,g}\right)_{g\in G}$ are uncorrelated local square integrable martingales, 
    \begin{align*}
        \partial\left[Z_g,Z_h\right](t)
        &= \sum_{\ell \in G} H_{g,\ell}(t) H_{h,\ell}(t) \partial\left[ \tld{M}_{E,\ell}\right](t)\\
        &= \sum_{\ell \in G} H_{g,\ell}(t) H_{h,\ell}(t) \sum_{i \in \ell} \frac{\partial\left[ M_i\right](t)}{a_i(t)^2} \\
        &= \sum_{\ell \in G} H_{g,\ell}(t) H_{h,\ell}(t) \sum_{i \in \ell} \frac{\partial N_i(t)}{a_i(t)^2}\\
        &= \partial \Gamma_{g,h}(t).
    \end{align*}

    \emph{Proof of \cref{itm:cv_Gamma}:} We want to prove that, $\forall g,h \in G$, $$\frac{\Gamma_{g,h}(t)}{n} \xrightarrow[n \to \infty]{\mathbb P} V_{g,h}(t).$$
    Recall that we just showed that $$\Gamma_{g,h}(t) = \sum_{\ell \in G} \int_{0}^t H_{g,\ell}(s) H_{h,\ell}(s) \sum_{i \in \ell} \frac{\partial N_i(s)}{a_i(s)^2}.$$

    Remember that our patients $(E_i,P_i,\bm X_i)$ are $i.i.d$, and assume (without loss of generality) that their group assignment are included in (or driven by) their covariates $\bm X$.
    Formally, consider that $G$ is a partition of the covariates support $\mathrm{Supp}(\bm X)$, and introduce $\forall \bm x \in \mathrm{Supp}(\bm X)$ the following notations to highlight this:
    \begin{align*}
        a_{\bm x}(s) &= \P{P \ge s}{E = s, \bm X = \bm x}\\
        b_{\bm x}(s) &= \P{P = s}{E \ge s, \bm X = \bm x}\\
        c_{\bm x}(s) &= \P{P \ge s}{E \ge s, \bm X = \bm x}\\
        N_{\bm x}(s) &= \I{O \le s, X \ge s \;|\; \bm X = \bm x}
    \end{align*}

    Let us first study the convergence of $R_g(t) = \frac{\tld{Y}_{E,g}(t)}{\tld{Y}_{E,\bigcdot}(t)}$.
    Note that : 
    \begin{align*}
        \Etm{\tld{Y}_{E,g}(t)}
            &= \frac{1}{n}\sum_{i=1}^n \Etm{\tld{Y}_{E,i}(t) \I{\bm X_i \in g}}\\
            &= \frac{1}{n}\sum_{i=1}^n Y_{E,i}(t) \I{\bm X_i \in g}\\
            &= \frac{1}{n}\sum_{i=1}^n \I{E_i \ge t, C_i \ge t, \bm X_i \in g}\\
            &\xrightarrow[n \to \infty]{\mathbb P} \P{E \ge t, C\ge t, \bm X \in g},
    \end{align*}
    and, likewise, $\Etm{\tld{Y}_{E,\bigcdot}(t)} \xrightarrow[n \to \infty]{\mathbb P} \P{E \ge t, C\ge t}$.


    Furthermore,
    \begin{align*}
        \V{\tld{Y}_{E,i}(t)}
            &= \V{\frac{Y_i(t)}{c_i(t)}}\\
            &= \frac{1}{c_i(t)^2} \V{Y_i(t)}\\
            &\le \frac{1}{c_i(t)^2}
    \end{align*}
    where the bounding is drastic but enough for our purposes.
    Indeed, we assumed that $c_{i}(t) \ge \epsilon$ over $t \in [0,T]$.
    We therefore have (by squeezing):
    $$\V{\tld{Y}_{E,\bigcdot}(t)}\xrightarrow[n \to \infty]{\mathbb P} 0,$$ 
    and similarly $\forall g \in G$, $$\V{\tld{Y}_{E,g}(t)}\xrightarrow[n \to \infty]{\mathbb P} 0.$$
    Hence: 
    \begin{align*}
        R_g(t) 
            \xrightarrow[n \to \infty]{\mathbb P} \frac{\P{E \ge t, C\ge t, \bm X \in g}}{\P{E \ge t, C\ge t}}
            = \P{\bm X \in g},
    \end{align*}
    where the last equality holds since censorship is independent and because we work under $\left(H_0\right)$.

    Finally, under $\left(H_0\right)$,
    $$H_{g,\{i\}}(s) = \I{\bm X_i \in g} - R_g(s) \xrightarrow[n \to \infty]{\mathbb P} \I{\bm X_i \in g} - \P{\bm X \in g},$$

    Denote now $$\Omega_{g,h}(\bm X) = \left(\I{\bm X \in g} - \P{\bm X \in g}\right)\left(\I{\bm X \in h} - \P{\bm X \in h}\right).$$ 

    Since $\left\lvert H_{g,\{i\}}(s) \right\rvert \le 1$, the dominated convergence theorem applies, and we have 
    \begin{align*}
        \lim_{n \to \infty} \int_{0}^t H_{g,\{i\}}(s) H_{h,\{i\}}(s) \frac{\partial N_i(s)}{a_i(s)^2}
        &= \int_{0}^t \Omega_{g,h}(\bm X_i) \frac{\partial N_i(s)}{a_i(s)^2}\\
    \end{align*}

    Thus,
    \begin{align*}
        \lim_{n \to\infty} \frac{\Gamma_{g,h}(t)}{n}
        &= \lim_{n \to\infty} \frac{1}{n} \sum_{\ell \in G} \int_{0}^t H_{g,\ell}(s) H_{h,\ell}(s) \sum_{i \in \ell} \frac{\partial N_i(s)}{a_i(s)^2}\\
        &= \lim_{n \to\infty} \frac{1}{n} \sum_{i=1}^{n} \int_{0}^t H_{g,\{i\}}(s) H_{h,\{i\}}(s) \frac{\partial N_i(s)}{a_i(s)^2}\\
        &= \E{\int_{0}^t \Omega_{g,h}(\bm X) \frac{\partial N_{\bm X}(s)}{a_{\bm X}(s)^2}}\\
        &= \E{\Omega_{g,h}(\bm X) \int_{0}^t \frac{\partial N_{\bm X}(s)}{a_{\bm X}(s)^2}}.
    \end{align*}

    We cannot unfortunately give a simpler expression for this limit, but we can however prove that it exists since we assumed $a_i(t) > \epsilon$ over $t \in [0,T]$.
    Since furthermore $\left\lvert \Omega_{g,\bm X}(s)\right\rvert \le 1$, and since $N_{\bm X}(t) \in {0,1}$, the expectation is positive and bounded above by $$\epsilon^{-2} \E{\int_{0}^t \partial N_{\bm X}(s)} = \epsilon^{-2} \E{N_{\bm X}(t)} \le \epsilon^{-2} < \infty.$$

    Call now $$V_{g,h}(t) = \E{\Omega_{g,h}(\bm X) \int_{0}^t \frac{\partial N_{\bm X}(s)}{a_{\bm X}(s)^2}}.$$

    Since $\bm \Gamma(t) = \left[\bm Z\right](t)$ is by construction semi-definite positive, $\bm V(t)$ also is  and \cref{itm:cv_Gamma} holds.

    \emph{Proof of \cref{itm:cv_Z}:} The Gaussian limit of $n^{-\frac{1}{2}}\bm Z$ comes directly from Robolledo's martingale central limit theorem, see \cite[Theorem II.5.1]{andersenCountingProcessModels1984}.
    There are two main conditions for this result to hold: 
    \begin{itemize}
        \item $\forall t, \left[ \frac{\bm Z}{\sqrt{n}}\right](t) \xrightarrow[n \to \infty]{\mathbb P} \bm V(t)$, where $\bm V(t)$ is a deterministic function with semi-definite positive matrix values, 
        \item $\forall t, \forall g\in G, \forall \eta > 0, \left\langle K_{g,\eta} \right\rangle(t) \xrightarrow[n \to \infty]{\mathbb P} 0$, where $K_{g,\eta}$ is a martingale that contains all jumps of $\frac{Z_g}{\sqrt{n}}$ of size greater than $\eta$ \emph{(Lindeberg's condition).}
    \end{itemize}

    The first condition is exactly \cref{itm:cv_Gamma} we just proved.
    The remaining Lindeberg condition can be verified as follows.
    Let us first construct $K_{g,\eta}$ as the restriction of $\frac{Z_g}{\sqrt{n}}$ to jumps of size greater than $\eta$.
    For that, recall from \cref{eq:Z} that, under the null,

    \begin{align*}
        \frac{Z_g(t)}{\sqrt{n}} 
        &= \int_{0}^t \frac{1}{\sqrt{n}}\sum_{\ell \in G} H_{g,\ell}(s)\partial \tld{M}_{E,\ell}(s)\\
        &= \int_{0}^t \frac{1}{\sqrt{n}}\sum_{\ell \in G} H_{g,\ell}(s)\partial \sum_{i \in \ell} \frac{\partial M_i(s)}{a_i(s)}\\
        &= \int_{0}^t \sum_{i =1}^n \frac{H_{g,\ell}(s)}{a_i(s)\sqrt{n}} \partial M_i(s) 
    \end{align*}

    Then, 
    since $N_i$'s have jumps of size 1, we can define the local square integrable martingale satisfying the requirement $K_{g,\eta}(t)$ as follows: 
    $$K_{g,\eta}(t) = \sum_{i =1}^n \int_{0}^t \frac{H_{g,\ell}(s)}{a_i(s)\sqrt{n}} \I{\left\lvert \frac{H_{g,\ell}(s)}{a_i(s)\sqrt{n}} \right\rvert > \eta} \partial M_i(s).$$
    

    Now, remembering that $\partial \langle M_i\rangle(s) = a_i(s)^2 \partial \Lambda_E(s)$, the predictable variate process of $K_{g,\eta}$ can be developed as:
    \begin{align}
        \left\langle K_{g,\eta}\right\rangle(t) 
        &= \sum_{i =1}^n \int_{0}^t \left(\frac{H_{g,\ell}(s)}{a_i(s)\sqrt{n}}\right)^2 \I{\left\lvert \frac{H_{g,\ell}(s)}{a_i(s)\sqrt{n}} \right\rvert > \eta} \partial \langle M_i\rangle(s)\nonumber\\
        &= \frac{1}{n}\sum_{i =1}^n  \int_{0}^t H_{g,\ell}(s)^2 \I{\left\lvert H_{g,\ell}(s) \right\rvert > \eta \sqrt{n} a_i(s)}\partial \Lambda_E(s)\nonumber\\
        &\le \frac{1}{n}\sum_{i =1}^n \int_{0}^t \I{\left\lvert H_{g,\ell}(s) \right\rvert > \eta \sqrt{n}  a_i(s)} \partial \Lambda_E(s)\label{eq:uperbound}\\
        &\xrightarrow[n \to \infty]{\mathbb P} 0,\label{eq:cvtozero}
    \end{align}
    where the upper bound in \cref{eq:uperbound} is holding since $\lvert H_{g,\ell}(s) \rvert \le 1$.
    Moreover, as soon as $n \ge \eta^{-2} a_i(s)^{-2} \ge \eta^{-2}\epsilon^{-2}$ (since we assumed $a_i(s) \ge \epsilon$ over $s \in [0,T]$), the indicatrix $$\I{\lvert H_{g,\ell}(s)\rvert > \epsilon\sqrt{n}a_i(s)}$$ is forced to be zero, and so the whole $\left\langle K_{g,\epsilon}\right\rangle(t)$ become zero after this particular $n$: the limit from \cref{eq:cvtozero} therefore holds and the Lindeberg's condition is fulfilled.

    \emph{Proof of \cref{itm:rank_V}:} We want to show that $\mathrm{Ker}(\bm V(t)) = \mathrm{Vect}(\bm 1)$.
    Recall the expressions for elements of $\bm V(t)$:
    $$V_{g,h}(t) = \E{\Omega_{g,h}(\bm X) \int_{0}^t \frac{\partial N_{\bm X}(s)}{a_{\bm X}(s)^2}}.$$
    Remark that the matrix $\bm \Omega(\bm X)$ can only take $\lvert G \rvert$ different values, depending on the group assignment of $\bm X$.
    Consider then the vectorized notations: 
    \begin{align*}
        g(\bm X) &= g \in G\text{ s.t. }\bm X \in g\\
        \bm i_g &= \left(\I{h =g}, h \in G\right)\\
        \bm p &=  \left(\P{\bm X \in g}, g \in G\right)\\
        \bm \Omega_g &= \left(\bm i_g - \bm p\right)\left(\bm i_g - \bm p\right)',
    \end{align*}
    so that we can write the matrix $\bm \Omega(\bm X)$ as $\bm \Omega_{g(\bm X)}$.
    Now, $\forall \bm u \in \mathbb R^{\lvert G \rvert}$,
    \begin{align*}
        \bm u \in \mathrm{Ker}(\bm \Omega_g)
            &\iff \bm u'\Omega_g\bm u = 0\\
            &\iff \bm u' \left(\bm i_g - \bm p\right)\left(\bm i_g - \bm p\right)' \bm u = 0\\
            &\iff \lVert \bm u' \left(\bm i_g - \bm p\right) \rVert_2^2 = 0\\
            &\iff u_g = \bm u'\bm p
    \end{align*}
    Moreover, we can decompose $\bm V(t)$ as: 
    \begin{align*}
        \bm V(t) 
            &= \E{\Omega(\bm X) \int_{0}^t \frac{\partial N_{ \bm X}(s)}{a_{\bm X}(s)^2}}\\
            &= \sum_{g \in G} \bm \Omega_g \E{\int_{0}^t \frac{\partial N_{ \bm X}(s)}{a_{\bm X}(s)^2}}{\bm X \in g} \P{\bm X \in g}\\
            &= \sum_{g \in G} w_g \bm \Omega_g,
    \end{align*}
    Assuming without loss of generality that the probability of uncensored deaths occurring before time $t$ is strictly positive in each group, the weights $w_g$ are all strictly positives. 
    Since moreover matrices $\Omega_g$ are all semi-definite positives,
    \begin{align*}
        \mathrm{Ker}(\bm V(t))
            &= \bigcap_{g \in G} \mathrm{Ker}(\bm \Omega_g)\\
            &= \bigcap_{g \in G} \{\bm u \in \mathbb R^{\lvert G \rvert}: \; u_g = \bm u'\bm p\}\\
            &= \{\bm u \in \mathbb R^{\lvert G \rvert}: \forall g \in G \; u_g = \bm u'\bm p\}\\
            &= \mathrm{Vect}(\bm 1).
    \end{align*}
\end{proofE}

\begin{propositionE}[Assymptotic test][end,restate]\label{prop:test} 
    Under $(H_0)$, assuming the existence of an $\epsilon > 0$ such that $a_i(t) > \epsilon$ and $c_i(t)>\epsilon$ over $t \in [0,T]$, we have
    $$\tld{\chi}(T) \xrightarrow[n \to\infty]{\mathcal D} \textup{\texttt{Chi2}}\left(\lvert G \rvert -1\right).$$ 
\end{propositionE}
\begin{proofE} 
    This result is directly drawn from \cref{lem:test}, from which we fetched the condition on the coefficients, since \cref{lem:test} gives us that $$n^{-\frac{1}{2}}\bm Z(t) \xrightarrow[n\to\infty]{\mathcal D} \N{0,\bm V(t)}.$$  From \cref{lem:test} again, we have $\mathrm{rank}(\bm V(t)) = \lvert G \rvert - 1$, the number of degree of freedom in the resulting $\chi^2$ distribution is $\lvert G \rvert - 1$.
\end{proofE}
The statistic $\tld{\chi}(T)$ can therefore be used to reject $\left(H_0\right)$ when $\tld{\chi}(T)$ gets above the $\alpha$-quantile of the $\texttt{Chi2}(\lvert G \rvert - 1)$ distribution, and such a test will have an asymptotic level $\alpha$.

As previously, $\tld{\chi}(T)$ is observable under $\HPi$ but not under $\HC$, and thus this test can only be used to possibly reject $H_0$ under $\HPi$.
However, this time, the plug-in to obtain an observable test statistic under $\HC$ must be carefully expressed.

Denote by $\hat{Y}_{E,g}, \hat{N}_{E,g}, \hat{\sigma}_{E,g}^2$ the estimators from \cref{eq:gen_ppe,eq:gen_var} computed on only observations from the group $g$.
Remark that the stochastic processes $\hat{a}_i, \hat{b}_i, \hat{c}_i$ will not be the same depending on the group.
In particular, the global estimates will \emph{not} match the sum of groupwise ones:
\begin{align*}
    \partial\hat{N}_{E}(t) &\neq \partial\hat{N}_{E,\bigcdot}(t) = \sum_{g \in G} \partial\hat{N}_{E,g}(t)\\
    \hat{Y}_{E}(t) &\neq \hat{Y}_{E,\bigcdot}(t) =  \sum_{g \in G} \hat{Y}_{E,g}(t).
\end{align*}

\begin{definition}[Observable test statistic] We propose to use the test statistic $\hat{\chi}(T) = \hat{\bm Z}(T)'\hat{\bm \Gamma}(T)^{-1} \hat{\bm Z}(T)$, where
    \begin{align*}
        \hat{R}_g(t) &= \frac{\hat{Y}_{E,g}(t)}{\hat{Y}_{E,\bigcdot}(t)}\\
        \hat{Z}_{g}(t) &= \hat{N}_{E,g}(t) - \int_{0}^t \hat{R}_g(s) \partial \hat{N}_{E,\bigcdot}(s)\\
        \hat{\Gamma}_{g,h}(t) &= \sum_{\ell\in G} \int_{0}^t \left(\delta_{\ell,g} - \hat{R}_g(s)\right)\left(\delta_{\ell,h} - \hat{R}_h(s)\right) \hat{Y}_{E,\ell}^2(s) \hat{\sigma}_{E,\ell}^2(s).
    \end{align*}
\end{definition}

In \cref{sec:simus,sec:examples}, we try to assess the performance of the excess hazard estimator, the estimator of its variance and last but not least the test statistic on simulated data through several scenarios.

\section{Simulation studies}\label{sec:simus}

We implemented the developed estimators in the pre-existing \texttt{Julia} package \texttt{NetSurvival.jl}~\cite{alhajal2024netsurvival}, alongside the reference implementations under $\HPi$.
Please see the package's documentation and code for technical details.
Furthermore, reproducing code for the below results (figures and tables from the simulations and examples sections) is available on a public GitHub repository\footnote{This repository: \url{https://www.github.com/lrnv/code_paper_netsurv_dep}}.

In simulations studies, we fix ourselves (and thus have access to) the true distribution of $E$ and the true dependence structure of $(E,P)$, described by the survival copula $\mathcal C_0$.
Furthermore, while sampling the model, we do observe for each individual the random variables $E_i,P_i,C_i$ and we can order them.
We leverage here this controlled environment to assess the performance and properties of our estimators.

\paragraph{Distributions of $\bm X$ and $P\;|\;\bm X$} We sample the covariates associated to our simulated patients from a very simple model: their sex is uniformly sampled between male and female, their diagnosis date uniformly between $1990$ and $2010$, and the age at diagnosis uniformly from $35$ to $75$ years old.
These demographic covariates determine the distribution of $P_i$'s through the Slovene rate table \texttt{slopop}\footnote{\texttt{RateTable} objects are implemented in the \texttt{RateTables.jl} Julia package.
They essentially map covariates (sex, date of diagnosis and age at diagnosis, sometimes geographical covariates too) to the distribution of $P$.
Technical details are extensively discussed in \cite{alhajal2024netsurvival}.
The \texttt{slopop} rate table is sourced from the R package \texttt{relsurv}~\cite{PermePavlik2018}.}.

\paragraph{Distribution of $E$} The first experiment considers a very simple true excess mortality following a $\texttt{Exponential}(\mu = 10)$ distribution (i.e., a constant excess hazard rate), for every individual.

\paragraph{Censorship $C$} We apply a censoring by $C \sim \texttt{Exponential}(\mu = 20) \wedge 15$ on top of this generation process, i.e., a random censoring and an administrative censoring at 15 years. This should yield a censoring rate of about one third. 

\begin{figure}[H]
    \centering\includegraphics[width=\textwidth]{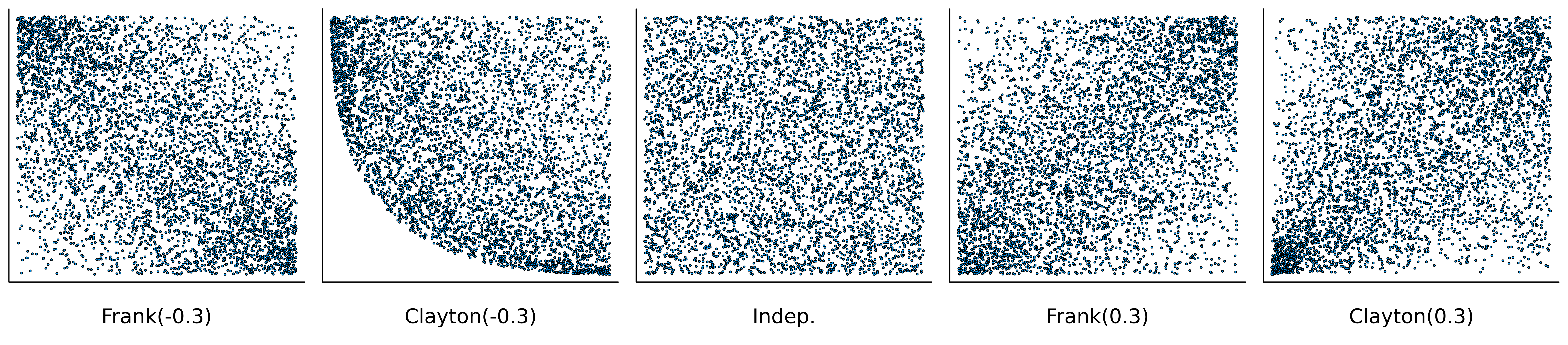}
    \caption{\label{fig:cops/showoff_cops} $5000$ pairs of ranks sampled (in the unit square) from different Archimedean copulas.
    We use these copulas as survival copulas of the vector $(E,P)$: the lower-left corner of each plot therefore shows the density of very large times $E$ and $P$, while the upper-right corner represents the density of very small times $E$ and $P$.
    Lower-right and upper-left corners respectively represent density of cases where $E$ is small while $P$ is large and vice versa.}
\end{figure}

\paragraph{True copula $\mathcal C_0$} The true survival copula $\mathcal C_0$ generating the data is taken among five copulas from the Archimedean family: $\texttt{Frank}(\tau = -0.3)$, $\texttt{Frank}(\tau = 0.3)$, $\texttt{Clayton}(\tau = -0.3)$, $\texttt{Clayton}(\tau = 0.3)$, and the independence copula $\Pi$.
Proper definition of these copulas is given in \cref{apx:distributions}.
A $5000$-sample of each of these copulas is displayed in \cref{fig:cops/showoff_cops}.

The Frank and Clayton families are parametrized by their Kendall $\tau$'s (see \cref{apx:distributions}), measuring the \emph{strength} of dependency on a $[-1;1]$ scale.
Remark on \cref{fig:cops/showoff_cops} the very different behavior of positive and negative parametrization in each family.
Even if the overall strength of dependency is the same, Frank and Clayton copulas display very different asymptotics in all corners, which correspond to different tail dependencies.
For example, the negative dependent Clayton copula has a support strictly smaller than the unit square and a frontier appears, implying (when used as a survival copula for $(E,P)$) that $E$ and $P$ will never be jointly large.
The strength of the dependency of our examples ($\tau \pm 0.3$) was chosen arbitrarily: see the next sections for a bit more insights about this parameter.

Let us take the example of the effect of age at diagnosis on breast cancer patients. According to \cite{adami1986relation}, or more recently \cite{eerola2001survival}, younger patients (with potentially higher values of $P$) have worse relative survival (that is lower values of $E$). Without more covariates to explain such a behavior, the Kendall $\tau$ of the dependence structure would probably be negative.

\paragraph{Hypothesized copula $\mathcal C$} We will perform estimation of the excess survival function and of its variance under $\HC$ for several potential dependence specifications $\mathcal C$.
In particular, we want to verify the amplitude of a potential wrong specification bias due to $\mathcal C \neq \mathcal C_0$ in the derivation of the estimator from \cref{eq:gen_ppe}.
We therefore pick $\mathcal C$, as $\mathcal C_0$, among the $5$ copulas from the upper list, we thus have $25$ different couples $(\mathcal C_0, \mathcal C)$.

\begin{figure}[H]
    \centering\includegraphics[width=\textwidth]{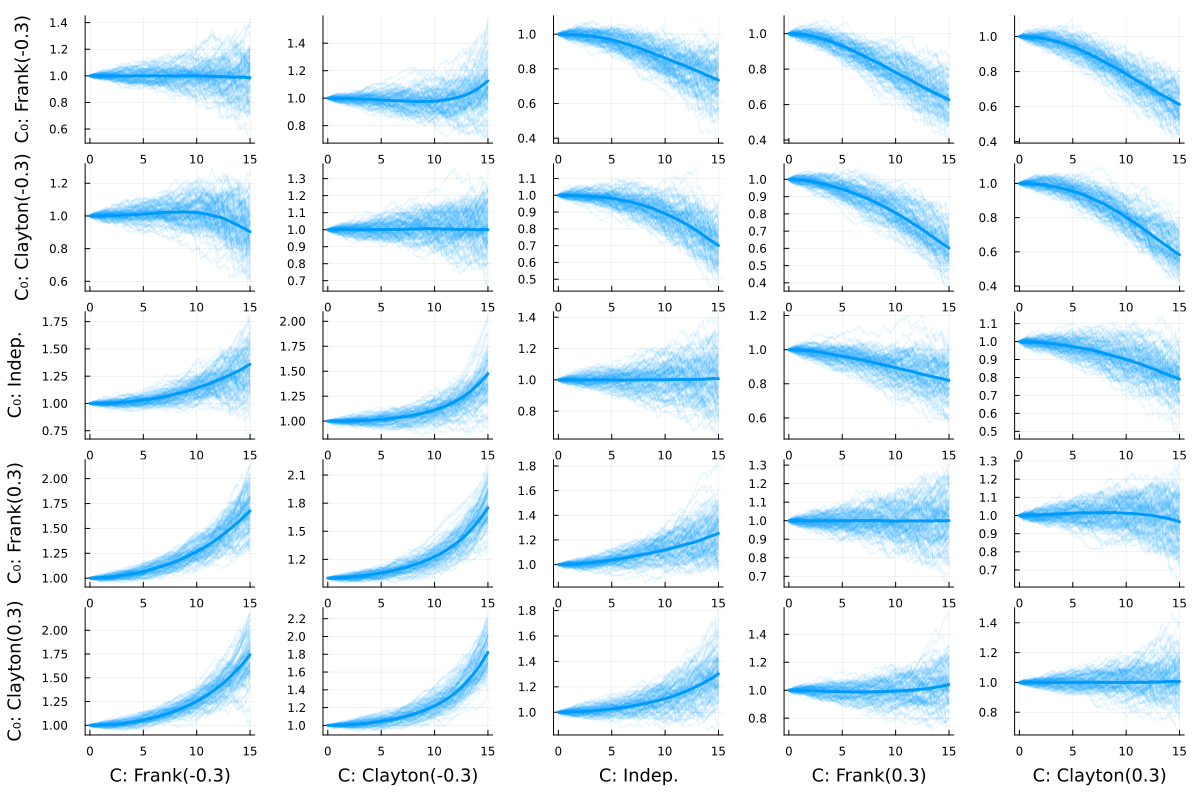}
    \caption{\label{fig:simus/featherplot} Each graph shows the $N=1000$ ratios $t \to \frac{\hat{S}_E^{(k)}(t)}{S_E(t)},\; k \in 1,...,N$ between the estimated survival curves and the true one.
    Abscissa represents time $t$ in years.
    Lines represent the true copula $\mathcal C_0$ and columns the hypothesized one $\mathcal C$: diagonal plots thus correspond to well-specified cases $\mathcal C = \mathcal C_0$.
    A flat average line close to 1 denotes a well-performing estimator.}
\end{figure}

For each of the potential $(\mathcal C_0,\mathcal C)$ case, we sample $N=1000$ datasets of size $n=500$ patients each. For each dataset index $k \in 1,...,N$, we denote $\hat{S}_E^{(k)}(t)$ the obtained excess survival function (estimated according to \cref{eq:gen_ppe}) and $\hat{\sigma}_E^{(k)}$ its variance (estimated according to \cref{eq:gen_var}). 
The ratios $$\frac{\hat{S}_E^{(k)}(t)}{S_E(t)}$$ are depicted scenario by scenario in \cref{fig:simus/featherplot}.
Straight up values of $\hat{S}_E^{(k)}(t)$ are also plotted in \cref{apx:extra_material}.

\begin{textAtEnd}[category=ExtraMaterial, end]
\begin{figure}[H]
    \centering\includegraphics[width=\textwidth]{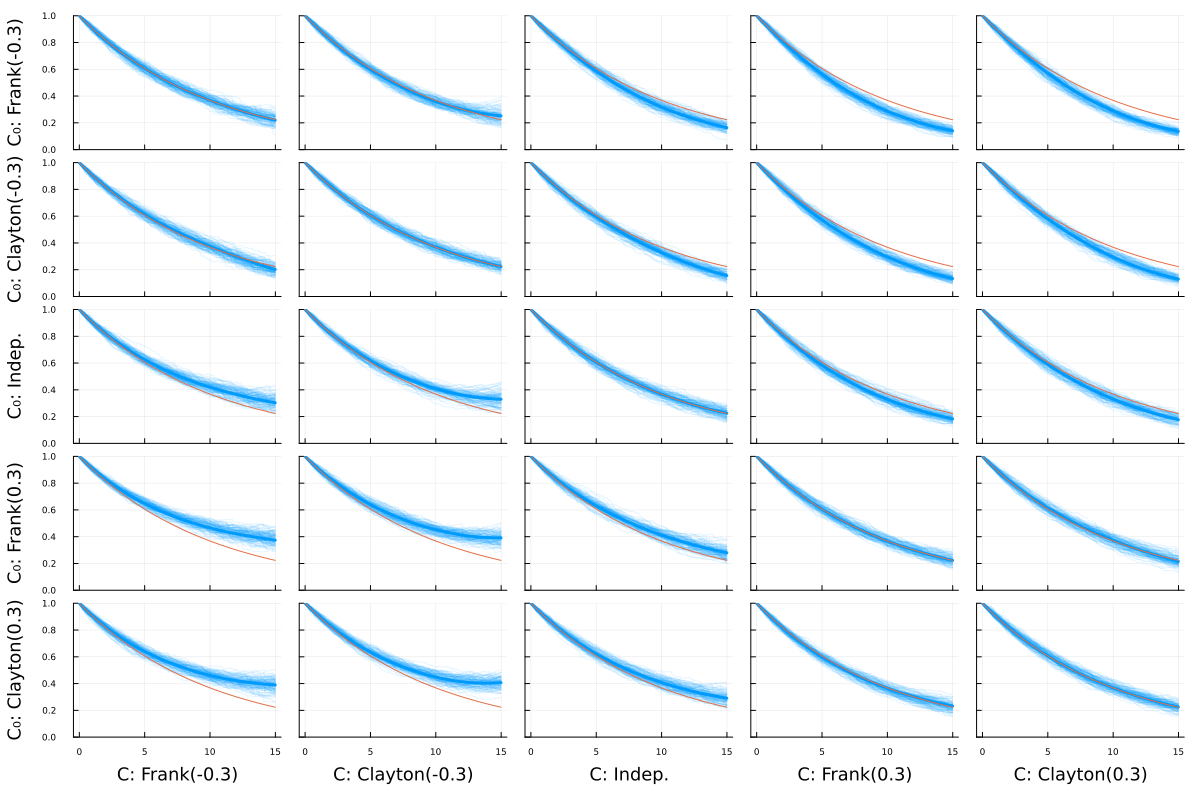}
    \caption{\label{fig:simus/featherplot2} Each graph shows (in blue) the $N=1000$ survival curves $t \to \hat{S}_E^{(1)}(t),...,t \to \hat{S}_E^{(N)}(t)$ estimated under $\HC$ on resamples simulated using $\mathcal C_0$.
    Abscissa represents time $t$ in years.
    The true survival curve $S_E(t)$ of $E \sim\texttt{Exponential}(\mu=10)$ is depicted in orange as a target reference.
    The true dependence structure $\mathcal C_0$ varies with the line, while the dependence structure assumed by the estimators $\mathcal C$ varies across columns: plots on the diagonal thus correspond to well-specified cases $\mathcal C = \mathcal C_0$.}
\end{figure}
\end{textAtEnd}

On the diagonal of \cref{fig:simus/featherplot} are well-specified cases where $\mathcal C = \mathcal C_0$.
First, remark that these well-specified cases yield (visually) unbiased estimation of the excess survival function, as we expected.
Second, remark that the bias direction looks to be driven by the difference of dependency strengths: the bias is positive when $\tau(\mathcal C_0) > \tau(\mathcal C)$ and negative when $\tau(\mathcal C_0) < \tau(\mathcal C)$.
In particular, the error obtained when using a wrong family with the right Kendall's $\tau$ is different from the error obtain when $\tau$ is wrong.

A more quantitative view of these results is given, for specific times $t \in [0,T]$, by the average bias, the root-mean-square error and the empirical coverage rate of the survival curves, defined by: 
\begin{align*}
    \mathrm{biais}(t) &= \frac{1}{N}\sum_{k=1}^N \hat{S_E}^{(k)}(t) - S_E(t)\\
    \mathrm{rmse}(t) &= \sqrt{\frac{1}{N}\sum_{k=1}^N \left(\hat{S_E}^{(k)}(t) - S_E(t)\right)^2}\\
    \mathrm{ecr}(t) &= \frac{1}{N} \sum_{k=1}^N \I{\left\lvert\ln S_E(t) - \ln \hat{S_E}^{(k)}(t)\right\rvert \le z_{0.975} \hat{\sigma_E}^{(k)}(t)}, 
\end{align*}
where $z_{\alpha}$ is the $\alpha$-quantile of a $\texttt{Normal}(0,1)$ distribution.
Values of $\mathrm{biais}(t)$, $\mathrm{rmse}(t)$ and $\mathrm{ecr}(t)$ can be observed at terminal time $t=15$ in \cref{tab:simus/rez15}.
Same extracts at earlier times $t\in \{5,10\}$ are available in \cref{apx:extra_material}.

\begin{textAtEnd}[category=ExtraMaterial, end]
    \begin{table}[H]
        \centering\footnotesize
        \caption{\label{tab:simus/rez5} Observed bias, root-mean-square error and empirical coverage rate over the $N=1000$ simulations at time $t=5$.
        Lines represent the true copula $\mathcal C_0$ and columns the hypothesized one $\mathcal C$.
        Smaller bias and root-mean-square errors are better, while empirical coverage rates closer to $95\%$ are better.}
        \begin{tabularx}{\linewidth}{lYYYYY}
  \toprule
   & $\texttt{Frank}(-0.3)$ & $\texttt{Clayton}(-0.3)$ & Indep. & $\texttt{Frank}(0.3)$ & $\texttt{Clayton}(0.3)$ \\\midrule
  $\mathrm{biais}(5)$ &  &  &  &  &  \\
  $\texttt{Frank}(-0.3)$ & $0.0005$ & $-0.0078$ & $-0.0199$ & $-0.0423$ & $-0.0354$ \\
  $\texttt{Clayton}(-0.3)$ & $0.0072$ & $0.0004$ & $-0.0108$ & $-0.0346$ & $-0.0275$ \\
  Indep. & $0.0186$ & $0.0110$ & $-0.0003$ & $-0.0235$ & $-0.0174$ \\
  $\texttt{Frank}(0.3)$ & $0.0410$ & $0.0335$ & $0.0225$ & $-0.0000$ & $0.0076$ \\
  $\texttt{Clayton}(0.3)$ & $0.0335$ & $0.0274$ & $0.0157$ & $-0.0068$ & $0.0004$ \\
  &&&&&\\[-5pt] $\mathrm{rmse}(5)$ &  &  &  &  &  \\
  $\texttt{Frank}(-0.3)$ & $0.0251$ & $0.0265$ & $0.0332$ & $0.0495$ & $0.0439$ \\
  $\texttt{Clayton}(-0.3)$ & $0.0263$ & $0.0248$ & $0.0281$ & $0.0434$ & $0.0377$ \\
  Indep. & $0.0307$ & $0.0271$ & $0.0254$ & $0.0350$ & $0.0321$ \\
  $\texttt{Frank}(0.3)$ & $0.0481$ & $0.0417$ & $0.0347$ & $0.0262$ & $0.0269$ \\
  $\texttt{Clayton}(0.3)$ & $0.0424$ & $0.0372$ & $0.0300$ & $0.0265$ & $0.0254$ \\
  &&&&&\\[-5pt] $\mathrm{ecr}(5)$ &  &  &  &  &  \\
  $\texttt{Frank}(-0.3)$ & $0.9610$ & $0.9470$ & $0.8940$ & $0.6490$ & $0.7640$ \\
  $\texttt{Clayton}(-0.3)$ & $0.9430$ & $0.9620$ & $0.9410$ & $0.7520$ & $0.8400$ \\
  Indep. & $0.8750$ & $0.9220$ & $0.9530$ & $0.8610$ & $0.8890$ \\
  $\texttt{Frank}(0.3)$ & $0.5980$ & $0.7230$ & $0.8220$ & $0.9400$ & $0.9310$ \\
  $\texttt{Clayton}(0.3)$ & $0.7030$ & $0.8010$ & $0.8940$ & $0.9470$ & $0.9550$ \\\bottomrule
\end{tabularx}

    \end{table}
    \begin{table}[H]
        \centering\footnotesize
        \caption{\label{tab:simus/rez10} Observed bias, root-mean-square error and empirical coverage rate over the $N=1000$ simulations at time $t=10$.
        Lines represent the true copula $\mathcal C_0$ and columns the hypothesized one $\mathcal C$.
        Smaller bias and root-mean-square errors are better, while empirical coverage rates closer to $95\%$ are better.}
        \begin{tabularx}{\linewidth}{lYYYYY}
  \toprule
   & $\texttt{Frank}(-0.3)$ & $\texttt{Clayton}(-0.3)$ & Indep. & $\texttt{Frank}(0.3)$ & $\texttt{Clayton}(0.3)$ \\\midrule
  $\mathrm{biais}(10)$ &  &  &  &  &  \\
  $\texttt{Frank}(-0.3)$ & $-0.0003$ & $-0.0076$ & $-0.0506$ & $-0.0804$ & $-0.0787$ \\
  $\texttt{Clayton}(-0.3)$ & $0.0075$ & $0.0022$ & $-0.0403$ & $-0.0716$ & $-0.0730$ \\
  Indep. & $0.0523$ & $0.0412$ & $0.0002$ & $-0.0388$ & $-0.0370$ \\
  $\texttt{Frank}(0.3)$ & $0.0983$ & $0.0840$ & $0.0438$ & $-0.0008$ & $0.0044$ \\
  $\texttt{Clayton}(0.3)$ & $0.0931$ & $0.0806$ & $0.0396$ & $-0.0023$ & $-0.0003$ \\
  &&&&&\\[-5pt] $\mathrm{rmse}(10)$ &  &  &  &  &  \\
  $\texttt{Frank}(-0.3)$ & $0.0286$ & $0.0281$ & $0.0580$ & $0.0843$ & $0.0827$ \\
  $\texttt{Clayton}(-0.3)$ & $0.0303$ & $0.0286$ & $0.0489$ & $0.0761$ & $0.0774$ \\
  Indep. & $0.0604$ & $0.0497$ & $0.0305$ & $0.0471$ & $0.0465$ \\
  $\texttt{Frank}(0.3)$ & $0.1030$ & $0.0890$ & $0.0538$ & $0.0285$ & $0.0302$ \\
  $\texttt{Clayton}(0.3)$ & $0.0982$ & $0.0857$ & $0.0503$ & $0.0280$ & $0.0283$ \\
  &&&&&\\[-5pt] $\mathrm{ecr}(10)$ &  &  &  &  &  \\
  $\texttt{Frank}(-0.3)$ & $0.9650$ & $0.9650$ & $0.6230$ & $0.1750$ & $0.1930$ \\
  $\texttt{Clayton}(-0.3)$ & $0.9460$ & $0.9570$ & $0.7560$ & $0.2890$ & $0.2770$ \\
  Indep. & $0.5680$ & $0.7070$ & $0.9280$ & $0.7350$ & $0.7730$ \\
  $\texttt{Frank}(0.3)$ & $0.0880$ & $0.1630$ & $0.6450$ & $0.9540$ & $0.9370$ \\
  $\texttt{Clayton}(0.3)$ & $0.1090$ & $0.1960$ & $0.6930$ & $0.9510$ & $0.9510$ \\\bottomrule
\end{tabularx}

    \end{table}
\end{textAtEnd}

\begin{table}[H]
    \centering\footnotesize
    \caption{\label{tab:simus/rez15} Observed bias, root-mean-square error and empirical coverage rate over the $N=1000$ simulations at time $t=15$.
    Lines represent the true copula $\mathcal C_0$ and columns the hypothesized one $\mathcal C$.
    Smaller bias and root-mean-square errors are better, while empirical coverage rates closer to $95\%$ are better.}
    \begin{tabularx}{\linewidth}{lYYYYY}
  \toprule
   & $\texttt{Frank}(-0.3)$ & $\texttt{Clayton}(-0.3)$ & Indep. & $\texttt{Frank}(0.3)$ & $\texttt{Clayton}(0.3)$ \\\midrule
  $\mathrm{biais}(15)$ &  &  &  &  &  \\
  $\texttt{Frank}(-0.3)$ & $-0.0033$ & $0.0282$ & $-0.0593$ & $-0.0831$ & $-0.0866$ \\
  $\texttt{Clayton}(-0.3)$ & $-0.0217$ & $-0.0000$ & $-0.0671$ & $-0.0890$ & $-0.0931$ \\
  Indep. & $0.0801$ & $0.1063$ & $0.0016$ & $-0.0403$ & $-0.0468$ \\
  $\texttt{Frank}(0.3)$ & $0.1509$ & $0.1678$ & $0.0565$ & $0.0002$ & $-0.0077$ \\
  $\texttt{Clayton}(0.3)$ & $0.1662$ & $0.1837$ & $0.0674$ & $0.0090$ & $0.0015$ \\
  &&&&&\\[-5pt] $\mathrm{rmse}(15)$ &  &  &  &  &  \\
  $\texttt{Frank}(-0.3)$ & $0.0368$ & $0.0524$ & $0.0647$ & $0.0860$ & $0.0891$ \\
  $\texttt{Clayton}(-0.3)$ & $0.0374$ & $0.0281$ & $0.0710$ & $0.0913$ & $0.0953$ \\
  Indep. & $0.0922$ & $0.1164$ & $0.0335$ & $0.0478$ & $0.0529$ \\
  $\texttt{Frank}(0.3)$ & $0.1575$ & $0.1734$ & $0.0678$ & $0.0292$ & $0.0316$ \\
  $\texttt{Clayton}(0.3)$ & $0.1732$ & $0.1892$ & $0.0772$ & $0.0322$ & $0.0299$ \\
  &&&&&\\[-5pt] $\mathrm{ecr}(15)$ &  &  &  &  &  \\
  $\texttt{Frank}(-0.3)$ & $0.9460$ & $0.8030$ & $0.5400$ & $0.0980$ & $0.0660$ \\
  $\texttt{Clayton}(-0.3)$ & $0.9670$ & $0.9760$ & $0.3890$ & $0.0500$ & $0.0370$ \\
  Indep. & $0.4000$ & $0.1760$ & $0.9400$ & $0.7110$ & $0.6450$ \\
  $\texttt{Frank}(0.3)$ & $0.0470$ & $0.0040$ & $0.5680$ & $0.9430$ & $0.9390$ \\
  $\texttt{Clayton}(0.3)$ & $0.0300$ & $0.0030$ & $0.4470$ & $0.9180$ & $0.9490$ \\\bottomrule
\end{tabularx}

\end{table}

Biases and root-mean-square errors in \cref{tab:simus/rez15} basically confirm the graphical analysis we had on \cref{fig:simus/featherplot}: well-specified estimators are unbiased, while badly-specified ones are biased.
Note that the bias direction seem to follow the Kendall's $\tau$ miss specification signs, an observation which will be corroborated by the next section but has to stay speculative for now on.

Moreover, the empirical coverage rates $\mathrm{ecr}(15)$ gives us a first look at the bias of the variance estimator $\hat{\sigma_E}^2(t)$, and in particular at the impact of our plug-in estimation:
\begin{itemize}
    \item On well-specified cases (when $\mathcal C = \mathcal C_0$), we observe empirical coverage rates close to the $95\%$ target (consistently across time, see \cref{tab:simus/rez5,tab:simus/rez10} in \cref{apx:extra_material}).
    In particular, the rate observed in these cases are similar to the Pohar Perme rate ($\mathcal C = \mathcal C_0 = \Pi$), which suggests that the error induced by the plug-in is not so important.
    \item On badly-specified cases (when $\mathcal C \neq \mathcal C_0$), we observe very low empirical coverage rates, clearly off-target (with an increasing bias over time, see \cref{tab:simus/rez5,tab:simus/rez10} in \cref{apx:extra_material}).
    This is coherent with the visual inspection of \cref{fig:simus/featherplot}, and clearly shows the importance of the dependency assumption (recall here that $\partial\hat{\Lambda}$'s value at $t$ depends on its path before $t$).
\end{itemize}

\paragraph{Log-rank tests benchmarks} This second experiment considers two scenarios designed to assess the performance of our generalized, dependency-aware, excess log-rank test.
Each scenario uses $N=1000$ resamples of datasets of $500$ patients sampled with the same demographic covariates, population mortality and underlying true dependence structure $\mathcal C_0$ as in the first experiment.
Each of these samples is split into two equally sized groups ($250$ patients each), differently in each scenario:
\begin{itemize}
    \item In Scenario 1, the joint distribution of $(P,\bm X)$ is the same in each group (and correspond to the settings of the previous experiment).
    \item In Scenario 2, the first group has age between $35$ and $65$, and the second between $65$ and $75$.
\end{itemize}
Finally, patients have excess mortality $\texttt{Exponential}(\mu)$ as previously, with a parameter $\mu$ that now depends on the groups: 
\begin{itemize}
    \item Under $H_0$, all $n = 500$ patients have $\mu=5$.
    \item Under $H_1$, the first group has $\mu=5$ while the second group has $\mu=6$ (hence a hazard ratio of $1.2$ between the two groups).
    \item Under $H_2$, the first group has $\mu=5$ while the second group has $\mu=10$ (hence a hazard ratio of $2$).
\end{itemize}

We therefore consider only proportional alternatives for the test, to stay in the scope of this paper.
Indeed, it is already known (see, e.g., \cite{wolskiPermutationTestBased2020}) that log-rank-type tests are not suitable under non-proportional alternatives, where other tests can perform better.
The proposition of a test allowing for non-proportional alternatives is left to further analysis.
The true dependence structure $\mathcal C_0$ could also have been varied from group to group, but to simplify the exposition of results, we kept the same $\mathcal C_0$ in the two groups.

We perform log-rank-type tests under $\HC$ for different $\mathcal C$.
Observed rejection rates in scenario 1 (at a nominal level $\alpha=5\%$) under each hypothesis, alongside asymptotic confidence intervals are given in \cref{tab:simus/test1/tst15} for tests at the terminal time $T=15$.
Confidence intervals for a rejection rate $\hat{r}$ computed as the mean of a set of rejection indicatrix $\bm i = \left(i_1,...,i_N\right)$ are obtained through a Gaussian asymptotics as $\hat{r} \pm z_{1-\alpha/2}\sqrt{\frac{\hat{r}(1-\hat{r})}{N}}$.
Similar results for $T\in \{5,10\}$ are available in \cref{apx:extra_material}.

\begin{textAtEnd}[category=ExtraMaterial, end]
    \begin{table}[H]
        \centering\footnotesize
        \caption{\label{tab:simus/test1/tst5} (Scenario 1) Rejection rate at $\alpha=5\%, T=5$, with their (asymptotically Gaussian) confidence intervals.
        Lines represent the true copula $\mathcal C_0$ and columns the hypothesized one $\mathcal C$, in each of the alternatives $H_0,H_1$ and $H_2$.}
        \begin{tabularx}{\linewidth}{lYYYYY}
  \toprule
   & $\texttt{Frank}(-0.3)$ & $\texttt{Clayton}(-0.3)$ & Indep. & $\texttt{Frank}(0.3)$ & $\texttt{Clayton}(0.3)$ \\\midrule
  $H_0$ &  &  &  &  &  \\
  $\texttt{Frank}(-0.3)$ & $4.1 \pm 1.2$ & $5.1 \pm 1.4$ & $6.7 \pm 1.6$ & $5.0 \pm 1.4$ & $5.1 \pm 1.4$ \\
  $\texttt{Clayton}(-0.3)$ & $4.3 \pm 1.3$ & $5.5 \pm 1.4$ & $5.4 \pm 1.4$ & $4.7 \pm 1.3$ & $5.1 \pm 1.4$ \\
  Indep. & $4.7 \pm 1.3$ & $4.7 \pm 1.3$ & $5.6 \pm 1.4$ & $5.8 \pm 1.4$ & $5.3 \pm 1.4$ \\
  $\texttt{Frank}(0.3)$ & $4.7 \pm 1.3$ & $5.5 \pm 1.4$ & $4.0 \pm 1.2$ & $6.2 \pm 1.5$ & $6.2 \pm 1.5$ \\
  $\texttt{Clayton}(0.3)$ & $4.1 \pm 1.2$ & $6.4 \pm 1.5$ & $5.1 \pm 1.4$ & $5.0 \pm 1.4$ & $5.9 \pm 1.5$ \\
  &&&&&\\[-5pt] $H_1$ &  &  &  &  &  \\
  $\texttt{Frank}(-0.3)$ & $28.6 \pm 2.8$ & $27.4 \pm 2.8$ & $30.2 \pm 2.8$ & $30.8 \pm 2.9$ & $30.3 \pm 2.8$ \\
  $\texttt{Clayton}(-0.3)$ & $29.7 \pm 2.8$ & $27.3 \pm 2.8$ & $30.6 \pm 2.9$ & $31.5 \pm 2.9$ & $33.9 \pm 2.9$ \\
  Indep. & $25.3 \pm 2.7$ & $25.5 \pm 2.7$ & $29.9 \pm 2.8$ & $26.2 \pm 2.7$ & $28.5 \pm 2.8$ \\
  $\texttt{Frank}(0.3)$ & $26.3 \pm 2.7$ & $28.2 \pm 2.8$ & $29.7 \pm 2.8$ & $28.9 \pm 2.8$ & $29.6 \pm 2.8$ \\
  $\texttt{Clayton}(0.3)$ & $27.9 \pm 2.8$ & $28.4 \pm 2.8$ & $29.5 \pm 2.8$ & $25.1 \pm 2.7$ & $31.7 \pm 2.9$ \\
  &&&&&\\[-5pt] $H_2$ &  &  &  &  &  \\
  $\texttt{Frank}(-0.3)$ & $99.8 \pm 0.3$ & $99.7 \pm 0.3$ & $99.9 \pm 0.2$ & $100.0 \pm 0.0$ & $100.0 \pm 0.0$ \\
  $\texttt{Clayton}(-0.3)$ & $100.0 \pm 0.0$ & $99.7 \pm 0.3$ & $99.8 \pm 0.3$ & $99.9 \pm 0.2$ & $100.0 \pm 0.0$ \\
  Indep. & $99.9 \pm 0.2$ & $99.8 \pm 0.3$ & $99.9 \pm 0.2$ & $99.8 \pm 0.3$ & $99.9 \pm 0.2$ \\
  $\texttt{Frank}(0.3)$ & $100.0 \pm 0.0$ & $99.8 \pm 0.3$ & $99.8 \pm 0.3$ & $100.0 \pm 0.0$ & $100.0 \pm 0.0$ \\
  $\texttt{Clayton}(0.3)$ & $99.6 \pm 0.4$ & $99.4 \pm 0.5$ & $100.0 \pm 0.0$ & $99.6 \pm 0.4$ & $99.8 \pm 0.3$ \\\bottomrule
\end{tabularx}

    \end{table}
    \begin{table}[H]
        \centering\footnotesize
        \caption{\label{tab:simus/test1/tst10} (Scenario 1) Rejection rate at $\alpha=5\%, T=10$, with their (asymptotically Gaussian) confidence intervals.
        Lines represent the true copula $\mathcal C_0$ and columns the hypothesized one $\mathcal C$, in each of the alternatives $H_0,H_1$ and $H_2$.}
        \begin{tabularx}{\linewidth}{lYYYYY}
  \toprule
   & $\texttt{Frank}(-0.3)$ & $\texttt{Clayton}(-0.3)$ & Indep. & $\texttt{Frank}(0.3)$ & $\texttt{Clayton}(0.3)$ \\\midrule
  $H_0$ &  &  &  &  &  \\
  $\texttt{Frank}(-0.3)$ & $4.1 \pm 1.2$ & $4.1 \pm 1.2$ & $5.9 \pm 1.5$ & $5.2 \pm 1.4$ & $4.7 \pm 1.3$ \\
  $\texttt{Clayton}(-0.3)$ & $4.3 \pm 1.3$ & $4.1 \pm 1.2$ & $3.8 \pm 1.2$ & $4.4 \pm 1.3$ & $4.5 \pm 1.3$ \\
  Indep. & $4.3 \pm 1.3$ & $3.1 \pm 1.1$ & $5.4 \pm 1.4$ & $5.1 \pm 1.4$ & $4.3 \pm 1.3$ \\
  $\texttt{Frank}(0.3)$ & $4.4 \pm 1.3$ & $5.1 \pm 1.4$ & $4.6 \pm 1.3$ & $5.1 \pm 1.4$ & $6.1 \pm 1.5$ \\
  $\texttt{Clayton}(0.3)$ & $4.9 \pm 1.3$ & $4.6 \pm 1.3$ & $5.9 \pm 1.5$ & $4.5 \pm 1.3$ & $6.3 \pm 1.5$ \\
  &&&&&\\[-5pt] $H_1$ &  &  &  &  &  \\
  $\texttt{Frank}(-0.3)$ & $33.4 \pm 2.9$ & $29.9 \pm 2.8$ & $35.3 \pm 3.0$ & $34.2 \pm 2.9$ & $35.4 \pm 3.0$ \\
  $\texttt{Clayton}(-0.3)$ & $38.9 \pm 3.0$ & $32.0 \pm 2.9$ & $38.5 \pm 3.0$ & $36.7 \pm 3.0$ & $40.0 \pm 3.0$ \\
  Indep. & $31.8 \pm 2.9$ & $30.5 \pm 2.9$ & $36.2 \pm 3.0$ & $33.5 \pm 2.9$ & $34.6 \pm 2.9$ \\
  $\texttt{Frank}(0.3)$ & $32.2 \pm 2.9$ & $30.1 \pm 2.8$ & $35.8 \pm 3.0$ & $36.4 \pm 3.0$ & $36.8 \pm 3.0$ \\
  $\texttt{Clayton}(0.3)$ & $32.5 \pm 2.9$ & $29.9 \pm 2.8$ & $34.1 \pm 2.9$ & $32.3 \pm 2.9$ & $38.3 \pm 3.0$ \\
  &&&&&\\[-5pt] $H_2$ &  &  &  &  &  \\
  $\texttt{Frank}(-0.3)$ & $100.0 \pm 0.0$ & $100.0 \pm 0.0$ & $100.0 \pm 0.0$ & $100.0 \pm 0.0$ & $100.0 \pm 0.0$ \\
  $\texttt{Clayton}(-0.3)$ & $100.0 \pm 0.0$ & $100.0 \pm 0.0$ & $100.0 \pm 0.0$ & $100.0 \pm 0.0$ & $100.0 \pm 0.0$ \\
  Indep. & $99.9 \pm 0.2$ & $99.9 \pm 0.2$ & $100.0 \pm 0.0$ & $100.0 \pm 0.0$ & $100.0 \pm 0.0$ \\
  $\texttt{Frank}(0.3)$ & $100.0 \pm 0.0$ & $99.9 \pm 0.2$ & $100.0 \pm 0.0$ & $100.0 \pm 0.0$ & $100.0 \pm 0.0$ \\
  $\texttt{Clayton}(0.3)$ & $99.9 \pm 0.2$ & $100.0 \pm 0.0$ & $100.0 \pm 0.0$ & $99.9 \pm 0.2$ & $100.0 \pm 0.0$ \\\bottomrule
\end{tabularx}

    \end{table}
\end{textAtEnd}
\begin{table}[H]
    \centering\footnotesize
    \caption{\label{tab:simus/test1/tst15} (Scenario 1) Rejection rate at $\alpha=5\%, T=15$, with their (asymptotically Gaussian) confidence intervals.
    Lines represent the true copula $\mathcal C_0$ and columns the hypothesized one $\mathcal C$, in each of the alternatives $H_0,H_1$ and $H_2$.}
    \begin{tabularx}{\linewidth}{lYYYYY}
  \toprule
   & $\texttt{Frank}(-0.3)$ & $\texttt{Clayton}(-0.3)$ & Indep. & $\texttt{Frank}(0.3)$ & $\texttt{Clayton}(0.3)$ \\\midrule
  $H_0$ &  &  &  &  &  \\
  $\texttt{Frank}(-0.3)$ & $4.7 \pm 1.3$ & $5.0 \pm 1.4$ & $4.6 \pm 1.3$ & $5.3 \pm 1.4$ & $4.8 \pm 1.3$ \\
  $\texttt{Clayton}(-0.3)$ & $3.9 \pm 1.2$ & $3.6 \pm 1.2$ & $4.0 \pm 1.2$ & $4.4 \pm 1.3$ & $4.4 \pm 1.3$ \\
  Indep. & $5.2 \pm 1.4$ & $7.7 \pm 1.7$ & $4.9 \pm 1.3$ & $5.1 \pm 1.4$ & $4.0 \pm 1.2$ \\
  $\texttt{Frank}(0.3)$ & $4.8 \pm 1.3$ & $9.2 \pm 1.8$ & $4.5 \pm 1.3$ & $5.5 \pm 1.4$ & $5.6 \pm 1.4$ \\
  $\texttt{Clayton}(0.3)$ & $7.1 \pm 1.6$ & $8.9 \pm 1.8$ & $7.1 \pm 1.6$ & $4.6 \pm 1.3$ & $5.6 \pm 1.4$ \\
  &&&&&\\[-5pt] $H_1$ &  &  &  &  &  \\
  $\texttt{Frank}(-0.3)$ & $32.6 \pm 2.9$ & $32.6 \pm 2.9$ & $36.2 \pm 3.0$ & $34.4 \pm 2.9$ & $36.7 \pm 3.0$ \\
  $\texttt{Clayton}(-0.3)$ & $40.8 \pm 3.0$ & $33.1 \pm 2.9$ & $39.2 \pm 3.0$ & $38.6 \pm 3.0$ & $40.8 \pm 3.0$ \\
  Indep. & $32.3 \pm 2.9$ & $30.2 \pm 2.8$ & $38.1 \pm 3.0$ & $35.8 \pm 3.0$ & $36.5 \pm 3.0$ \\
  $\texttt{Frank}(0.3)$ & $30.6 \pm 2.9$ & $30.5 \pm 2.9$ & $37.4 \pm 3.0$ & $35.9 \pm 3.0$ & $37.3 \pm 3.0$ \\
  $\texttt{Clayton}(0.3)$ & $32.9 \pm 2.9$ & $30.3 \pm 2.8$ & $35.5 \pm 3.0$ & $33.7 \pm 2.9$ & $40.5 \pm 3.0$ \\
  &&&&&\\[-5pt] $H_2$ &  &  &  &  &  \\
  $\texttt{Frank}(-0.3)$ & $100.0 \pm 0.0$ & $99.1 \pm 0.6$ & $100.0 \pm 0.0$ & $100.0 \pm 0.0$ & $100.0 \pm 0.0$ \\
  $\texttt{Clayton}(-0.3)$ & $100.0 \pm 0.0$ & $100.0 \pm 0.0$ & $100.0 \pm 0.0$ & $100.0 \pm 0.0$ & $100.0 \pm 0.0$ \\
  Indep. & $99.9 \pm 0.2$ & $98.0 \pm 0.9$ & $100.0 \pm 0.0$ & $100.0 \pm 0.0$ & $100.0 \pm 0.0$ \\
  $\texttt{Frank}(0.3)$ & $99.7 \pm 0.3$ & $98.8 \pm 0.7$ & $100.0 \pm 0.0$ & $100.0 \pm 0.0$ & $100.0 \pm 0.0$ \\
  $\texttt{Clayton}(0.3)$ & $99.8 \pm 0.3$ & $96.5 \pm 1.1$ & $100.0 \pm 0.0$ & $99.9 \pm 0.2$ & $100.0 \pm 0.0$ \\\bottomrule
\end{tabularx}

\end{table}

On \cref{tab:simus/test1/tst15}, rejection rates under $H_0$ represents the observed type 1 error of the test (target value is the nominal level $\alpha=0.05$), while the $H_1$ and $H_2$ rejection rates represent the test's observed power (bigger is better) against these alternatives.
As expected from e.g., \cite{graffeoLogRankType2016}, the power of the test is increasing with the hazard ratio and thus better results are obtained for $H_2$ than $H_1$, even when the dependence structure is wrongly specified.

\begin{figure}[H]
    \centering\includegraphics[width=\textwidth]{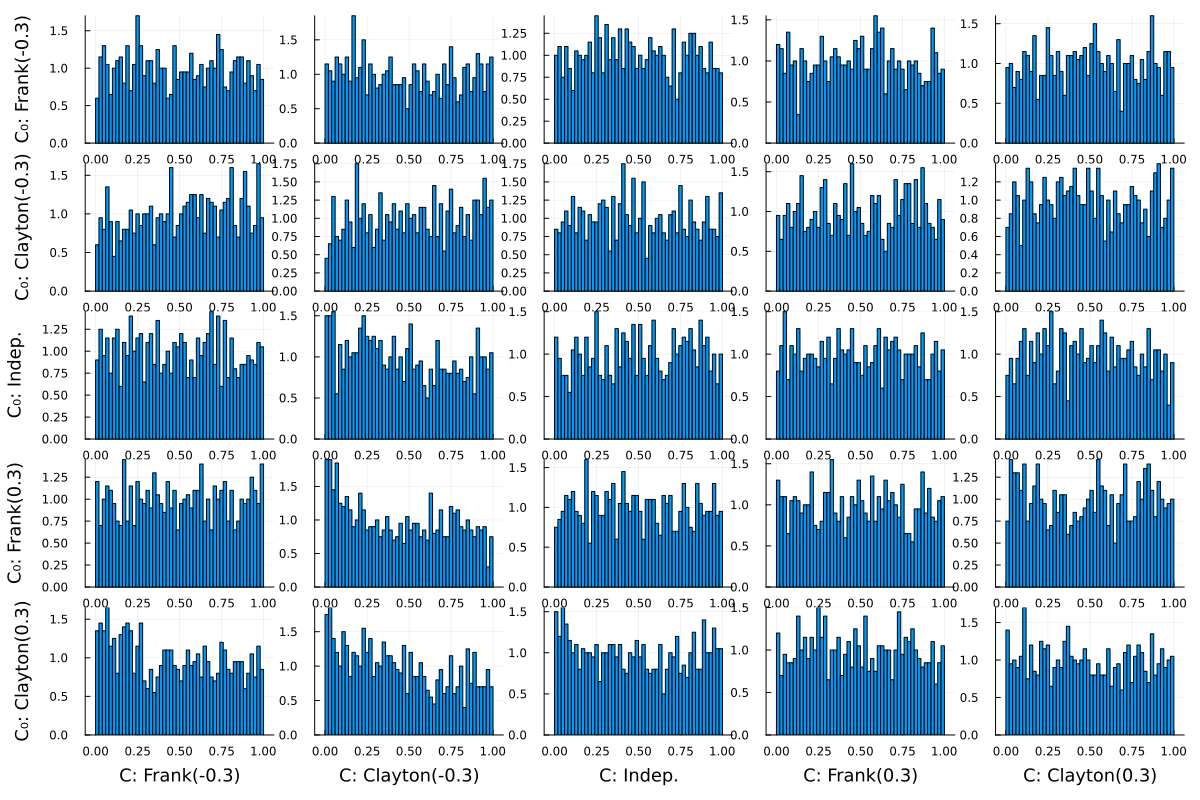}
    \caption{\label{fig:simus/test1/histo_H0} (Scenario 1) Each histogram represents the distribution of the $N=1000$ p-values under $H_0$.
    Lines represent the true copula $\mathcal C_0$ and columns the hypothesized one $\mathcal C$.
    A good results is a uniform distribution.}
\end{figure}

Moreover, in most of the $(\mathcal C_0, \mathcal C)$ specifications, the type 1 error is around the target $5\%$.
More complete information on the distributions of p-values in the $H_0$ case can be seen in \cref{fig:simus/test1/histo_H0}, where we plotted, in each $(\mathcal C_0, \mathcal C)$ case, a histogram of the p-values obtained at $T = 15$.

We first observe on \cref{fig:simus/test1/histo_H0} a slight issue when $\mathcal C_0$ is $\texttt{Frank}(0.3)$ or $\texttt{Clayton}(0.3)$ and estimating with $\mathcal C \in \left\{\texttt{Frank}(-0.3), \texttt{Clayton}(-0.3)\right\}$ (in the four lower-left histograms), but not really on other badly-specified cases.

This is coherent with the reading of \cref{tab:simus/test1/tst15}, and highlights a particular behavior of these estimators when $P$ is equally distributed in each group.
Indeed, Scenario 1 is favorable to bias in the copula specification, since the distribution of $P$ is each group is theoretically the same: even if the estimation of $E$ done in each group is highly biased, it has \emph{the same bias} in each group, hence yielding a somewhat acceptable type 1 error in \cref{tab:simus/test1/tst15} and somewhat acceptable p-value uniformity in \cref{fig:simus/test1/histo_H0}.

\begin{figure}[H]
    \centering\includegraphics[width=\textwidth]{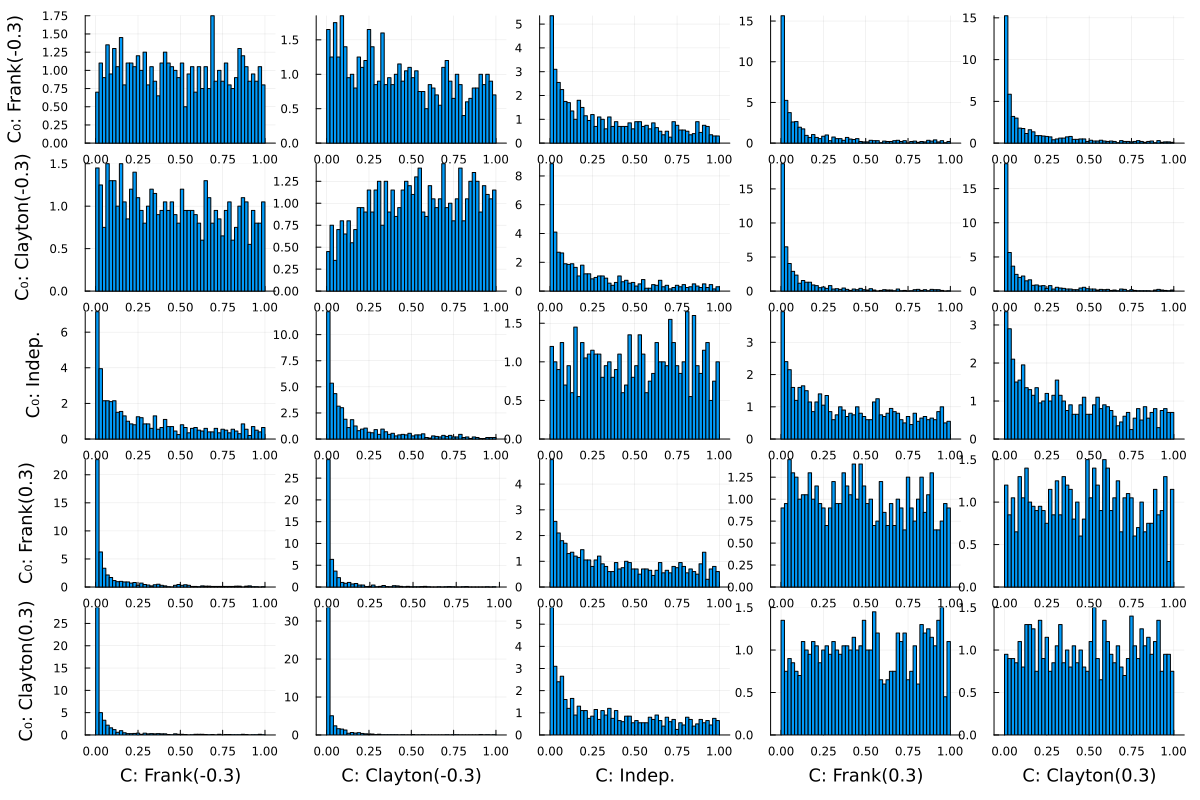}
    \caption{\label{fig:simus/test2/histo_H0} (Scenario 2) Each histogram represents the distribution of the $N=1000$ p-values under $H_0$.
    Lines represent the true copula $\mathcal C_0$ and columns the hypothesized one $\mathcal C$.
    A good results is a uniform distribution.}
\end{figure}

Let us now look at the same results for Scenario 2, where the distribution of $P$ is purposely not the same in each group (which is also a more realistic scenario).
For Scenario 2, \cref{tab:simus/test2/tst15} gives rejection rates at $T=15$, while rejection rates at $T\in \{5,10\}$ are available in \cref{apx:extra_material}, and \cref{fig:simus/test2/histo_H0} represents histograms of p-values under $H_0$.

\begin{textAtEnd}[category=ExtraMaterial, end]
    \begin{table}[H]
        \centering\footnotesize
        \caption{\label{tab:simus/test2/tst5} (Scenario 2) Rejection rate at $\alpha=5\%, T =5$, with their (asymptotically Gaussian) confidence intervals.
        Lines represent the true copula $\mathcal C_0$ and columns the hypothesized one $\mathcal C$, in each of the alternatives $H_0,H_1$ and $H_2$.}
        \begin{tabularx}{\linewidth}{lYYYYY}
  \toprule
   & $\texttt{Frank}(-0.3)$ & $\texttt{Clayton}(-0.3)$ & Indep. & $\texttt{Frank}(0.3)$ & $\texttt{Clayton}(0.3)$ \\\midrule
  $H_0$ &  &  &  &  &  \\
  $\texttt{Frank}(-0.3)$ & $4.9 \pm 1.3$ & $5.5 \pm 1.4$ & $10.8 \pm 1.9$ & $25.6 \pm 2.7$ & $22.5 \pm 2.6$ \\
  $\texttt{Clayton}(-0.3)$ & $4.4 \pm 1.3$ & $4.2 \pm 1.2$ & $9.1 \pm 1.8$ & $20.5 \pm 2.5$ & $18.0 \pm 2.4$ \\
  Indep. & $10.6 \pm 1.9$ & $6.2 \pm 1.5$ & $4.9 \pm 1.3$ & $10.6 \pm 1.9$ & $6.9 \pm 1.6$ \\
  $\texttt{Frank}(0.3)$ & $27.2 \pm 2.8$ & $20.1 \pm 2.5$ & $10.1 \pm 1.9$ & $4.4 \pm 1.3$ & $6.3 \pm 1.5$ \\
  $\texttt{Clayton}(0.3)$ & $22.8 \pm 2.6$ & $14.6 \pm 2.2$ & $9.7 \pm 1.8$ & $5.2 \pm 1.4$ & $4.8 \pm 1.3$ \\
  &&&&&\\[-5pt] $H_1$ &  &  &  &  &  \\
  $\texttt{Frank}(-0.3)$ & $29.5 \pm 2.8$ & $18.7 \pm 2.4$ & $10.3 \pm 1.9$ & $5.1 \pm 1.4$ & $5.0 \pm 1.4$ \\
  $\texttt{Clayton}(-0.3)$ & $36.4 \pm 3.0$ & $27.4 \pm 2.8$ & $15.1 \pm 2.2$ & $6.5 \pm 1.5$ & $7.3 \pm 1.6$ \\
  Indep. & $53.7 \pm 3.1$ & $45.8 \pm 3.1$ & $30.4 \pm 2.9$ & $10.8 \pm 1.9$ & $15.0 \pm 2.2$ \\
  $\texttt{Frank}(0.3)$ & $76.0 \pm 2.6$ & $68.8 \pm 2.9$ & $53.1 \pm 3.1$ & $27.2 \pm 2.8$ & $34.1 \pm 2.9$ \\
  $\texttt{Clayton}(0.3)$ & $72.9 \pm 2.8$ & $64.9 \pm 3.0$ & $47.4 \pm 3.1$ & $24.7 \pm 2.7$ & $29.4 \pm 2.8$ \\
  &&&&&\\[-5pt] $H_2$ &  &  &  &  &  \\
  $\texttt{Frank}(-0.3)$ & $99.9 \pm 0.2$ & $99.7 \pm 0.3$ & $98.4 \pm 0.8$ & $94.2 \pm 1.4$ & $96.0 \pm 1.2$ \\
  $\texttt{Clayton}(-0.3)$ & $100.0 \pm 0.0$ & $99.8 \pm 0.3$ & $99.1 \pm 0.6$ & $96.6 \pm 1.1$ & $97.7 \pm 0.9$ \\
  Indep. & $99.9 \pm 0.2$ & $100.0 \pm 0.0$ & $99.8 \pm 0.3$ & $98.1 \pm 0.8$ & $98.9 \pm 0.6$ \\
  $\texttt{Frank}(0.3)$ & $100.0 \pm 0.0$ & $100.0 \pm 0.0$ & $100.0 \pm 0.0$ & $100.0 \pm 0.0$ & $100.0 \pm 0.0$ \\
  $\texttt{Clayton}(0.3)$ & $100.0 \pm 0.0$ & $100.0 \pm 0.0$ & $99.9 \pm 0.2$ & $99.4 \pm 0.5$ & $99.8 \pm 0.3$ \\\bottomrule
\end{tabularx}

    \end{table}
    \begin{table}[H]
        \centering\footnotesize
        \caption{\label{tab:simus/test2/tst10} (Scenario 2) Rejection rate at $\alpha=5\%, T=10$, with their (asymptotically Gaussian) confidence intervals.
        Lines represent the true copula $\mathcal C_0$ and columns the hypothesized one $\mathcal C$, in each of the alternatives $H_0,H_1$ and $H_2$.}
        \begin{tabularx}{\linewidth}{lYYYYY}
  \toprule
   & $\texttt{Frank}(-0.3)$ & $\texttt{Clayton}(-0.3)$ & Indep. & $\texttt{Frank}(0.3)$ & $\texttt{Clayton}(0.3)$ \\\midrule
  $H_0$ &  &  &  &  &  \\
  $\texttt{Frank}(-0.3)$ & $5.2 \pm 1.4$ & $2.5 \pm 1.0$ & $18.4 \pm 2.4$ & $42.2 \pm 3.1$ & $41.3 \pm 3.1$ \\
  $\texttt{Clayton}(-0.3)$ & $4.5 \pm 1.3$ & $3.1 \pm 1.1$ & $23.2 \pm 2.6$ & $44.1 \pm 3.1$ & $45.1 \pm 3.1$ \\
  Indep. & $19.9 \pm 2.5$ & $17.8 \pm 2.4$ & $5.1 \pm 1.4$ & $13.2 \pm 2.1$ & $12.5 \pm 2.1$ \\
  $\texttt{Frank}(0.3)$ & $53.9 \pm 3.1$ & $49.6 \pm 3.1$ & $13.9 \pm 2.1$ & $5.2 \pm 1.4$ & $5.1 \pm 1.4$ \\
  $\texttt{Clayton}(0.3)$ & $55.5 \pm 3.1$ & $49.6 \pm 3.1$ & $15.0 \pm 2.2$ & $5.1 \pm 1.4$ & $4.4 \pm 1.3$ \\
  &&&&&\\[-5pt] $H_1$ &  &  &  &  &  \\
  $\texttt{Frank}(-0.3)$ & $32.7 \pm 2.9$ & $32.3 \pm 2.9$ & $7.0 \pm 1.6$ & $6.3 \pm 1.5$ & $4.8 \pm 1.3$ \\
  $\texttt{Clayton}(-0.3)$ & $29.1 \pm 2.8$ & $31.2 \pm 2.9$ & $7.4 \pm 1.6$ & $5.1 \pm 1.4$ & $5.7 \pm 1.4$ \\
  Indep. & $75.7 \pm 2.7$ & $73.9 \pm 2.7$ & $36.5 \pm 3.0$ & $11.7 \pm 2.0$ & $11.1 \pm 1.9$ \\
  $\texttt{Frank}(0.3)$ & $95.7 \pm 1.3$ & $94.5 \pm 1.4$ & $70.4 \pm 2.8$ & $34.6 \pm 2.9$ & $37.1 \pm 3.0$ \\
  $\texttt{Clayton}(0.3)$ & $95.3 \pm 1.3$ & $95.3 \pm 1.3$ & $69.9 \pm 2.8$ & $36.3 \pm 3.0$ & $36.0 \pm 3.0$ \\
  &&&&&\\[-5pt] $H_2$ &  &  &  &  &  \\
  $\texttt{Frank}(-0.3)$ & $100.0 \pm 0.0$ & $100.0 \pm 0.0$ & $98.7 \pm 0.7$ & $93.5 \pm 1.5$ & $92.9 \pm 1.6$ \\
  $\texttt{Clayton}(-0.3)$ & $99.9 \pm 0.2$ & $100.0 \pm 0.0$ & $99.7 \pm 0.3$ & $96.8 \pm 1.1$ & $97.1 \pm 1.0$ \\
  Indep. & $100.0 \pm 0.0$ & $100.0 \pm 0.0$ & $100.0 \pm 0.0$ & $99.4 \pm 0.5$ & $99.1 \pm 0.6$ \\
  $\texttt{Frank}(0.3)$ & $100.0 \pm 0.0$ & $100.0 \pm 0.0$ & $100.0 \pm 0.0$ & $100.0 \pm 0.0$ & $100.0 \pm 0.0$ \\
  $\texttt{Clayton}(0.3)$ & $100.0 \pm 0.0$ & $100.0 \pm 0.0$ & $100.0 \pm 0.0$ & $99.9 \pm 0.2$ & $99.9 \pm 0.2$ \\\bottomrule
\end{tabularx}

    \end{table}
\end{textAtEnd}
\begin{table}[H]
    \centering\footnotesize
    \caption{\label{tab:simus/test2/tst15} (Scenario 2) Rejection rate at $\alpha=5\%, T=15$, with their (asymptotically Gaussian) confidence intervals.
    Lines represent the true copula $\mathcal C_0$ and columns the hypothesized one $\mathcal C$, in each of the alternatives $H_0,H_1$ and $H_2$.}
    \begin{tabularx}{\linewidth}{lYYYYY}
  \toprule
   & $\texttt{Frank}(-0.3)$ & $\texttt{Clayton}(-0.3)$ & Indep. & $\texttt{Frank}(0.3)$ & $\texttt{Clayton}(0.3)$ \\\midrule
  $H_0$ &  &  &  &  &  \\
  $\texttt{Frank}(-0.3)$ & $4.6 \pm 1.3$ & $7.9 \pm 1.7$ & $19.8 \pm 2.5$ & $45.8 \pm 3.1$ & $45.6 \pm 3.1$ \\
  $\texttt{Clayton}(-0.3)$ & $5.9 \pm 1.5$ & $2.7 \pm 1.0$ & $28.8 \pm 2.8$ & $53.6 \pm 3.1$ & $51.9 \pm 3.1$ \\
  Indep. & $24.2 \pm 2.7$ & $39.1 \pm 3.0$ & $5.4 \pm 1.4$ & $15.0 \pm 2.2$ & $14.3 \pm 2.2$ \\
  $\texttt{Frank}(0.3)$ & $61.3 \pm 3.0$ & $75.8 \pm 2.7$ & $17.0 \pm 2.3$ & $5.3 \pm 1.4$ & $5.2 \pm 1.4$ \\
  $\texttt{Clayton}(0.3)$ & $70.2 \pm 2.8$ & $80.0 \pm 2.5$ & $19.7 \pm 2.5$ & $5.0 \pm 1.4$ & $4.6 \pm 1.3$ \\
  &&&&&\\[-5pt] $H_1$ &  &  &  &  &  \\
  $\texttt{Frank}(-0.3)$ & $33.3 \pm 2.9$ & $49.7 \pm 3.1$ & $5.8 \pm 1.4$ & $6.8 \pm 1.6$ & $6.6 \pm 1.5$ \\
  $\texttt{Clayton}(-0.3)$ & $13.8 \pm 2.1$ & $27.1 \pm 2.8$ & $5.0 \pm 1.4$ & $7.7 \pm 1.7$ & $6.6 \pm 1.5$ \\
  Indep. & $78.3 \pm 2.6$ & $89.6 \pm 1.9$ & $36.1 \pm 3.0$ & $10.3 \pm 1.9$ & $8.9 \pm 1.8$ \\
  $\texttt{Frank}(0.3)$ & $96.9 \pm 1.1$ & $98.9 \pm 0.6$ & $73.9 \pm 2.7$ & $37.5 \pm 3.0$ & $35.0 \pm 3.0$ \\
  $\texttt{Clayton}(0.3)$ & $97.9 \pm 0.9$ & $99.3 \pm 0.5$ & $78.6 \pm 2.5$ & $43.5 \pm 3.1$ & $39.0 \pm 3.0$ \\
  &&&&&\\[-5pt] $H_2$ &  &  &  &  &  \\
  $\texttt{Frank}(-0.3)$ & $99.9 \pm 0.2$ & $100.0 \pm 0.0$ & $98.4 \pm 0.8$ & $91.1 \pm 1.8$ & $90.4 \pm 1.8$ \\
  $\texttt{Clayton}(-0.3)$ & $99.9 \pm 0.2$ & $100.0 \pm 0.0$ & $99.2 \pm 0.6$ & $94.4 \pm 1.4$ & $94.0 \pm 1.5$ \\
  Indep. & $100.0 \pm 0.0$ & $100.0 \pm 0.0$ & $99.9 \pm 0.2$ & $99.4 \pm 0.5$ & $99.5 \pm 0.4$ \\
  $\texttt{Frank}(0.3)$ & $100.0 \pm 0.0$ & $100.0 \pm 0.0$ & $100.0 \pm 0.0$ & $100.0 \pm 0.0$ & $100.0 \pm 0.0$ \\
  $\texttt{Clayton}(0.3)$ & $100.0 \pm 0.0$ & $100.0 \pm 0.0$ & $100.0 \pm 0.0$ & $99.9 \pm 0.2$ & $99.9 \pm 0.2$ \\\bottomrule
\end{tabularx}

\end{table}

The $H_0$ and $H_1$ results on \cref{tab:simus/test2/tst15}, respectively corresponding to type 1 errors and powers of the test are striking in this second scenario.
When the distribution of $P$ is not the same in each group, we clearly observe (as expected) that proper dependence specifications $\mathcal C = \mathcal C_0$ gives correctly calibrated tests, while wrong specification induces a wrong type 1 error.
The distributions of p-values under $H_0$ in this scenario, in \cref{fig:simus/test2/histo_H0}, confirms this result.

The few off-diagonal histograms that are somewhat acceptable on \cref{fig:simus/test2/histo_H0} corresponds to cases where the dependence structure shape is wrongly specified, but its strength (measured by Kendall's $\tau$) is well-specified.
However, the results are better for $(\mathcal C_0,\mathcal C)$ being $\left(\texttt{Frank}(0.3),\texttt{Clayton}(0.3)\right)$ or 
$\left(\texttt{Clayton}(0.3),\texttt{Frank}(0.3)\right)$ than for $\left(\texttt{Frank}(-0.3),\texttt{Clayton}(-0.3)\right)$ or 
$\left(\texttt{Clayton}(-0.3),\texttt{Frank}(-0.3)\right)$.

\section{Applications}\label{sec:examples}

\begin{figure}[H]
    \centering
    \includegraphics[width=\textwidth]{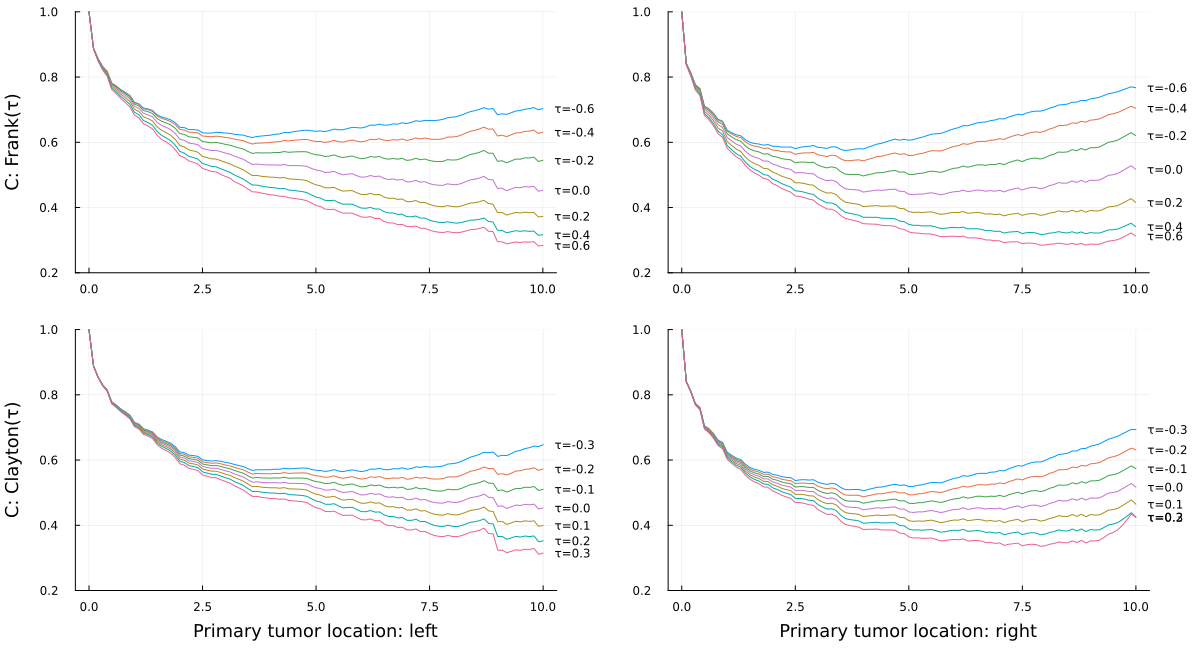}
    \caption{\label{fig:example/short_paper_graph} Estimated $t \to \hat{S}_E(t)$ for various $\HC$, time $t \in [0,T]$ is graduated in full years.
    Top graphs corresponds to Frank copulas while bottom are Clayton's, each curve corresponding to a different Kendall $\tau$.
    Left and right sides corresponds to the primary tumor location.}
\end{figure}

We use data on colorectal cancer patients extensively described in~\cite{giorgi2003relative,wolskiPermutationTestBased2020}.
This dataset contains one extra covariate for each patient that we will discuss here: the primary tumor location (left or right of the colon).
According to \cite{wolskiPermutationTestBased2020}, under $\HPi$, the left and right side net survival curves are crossing each other.
In such situations where hazard are non-proportionals, other tests can be more powerful, as they are able to detect differences that the log-rank type test cannot identify.
We try here to challenge the inconclusive result of the log-rank type test by looking at the dependence structure assumption.
Using several dependence structures, we first obtain net survival curves from Figure \ref{fig:example/short_paper_graph} on each of the two subgroups (according to primary tumor location).

\begin{textAtEnd}[category=ExtraMaterial, end]
    \begin{figure}[H]
        \centering
        \includegraphics[width=\textwidth]{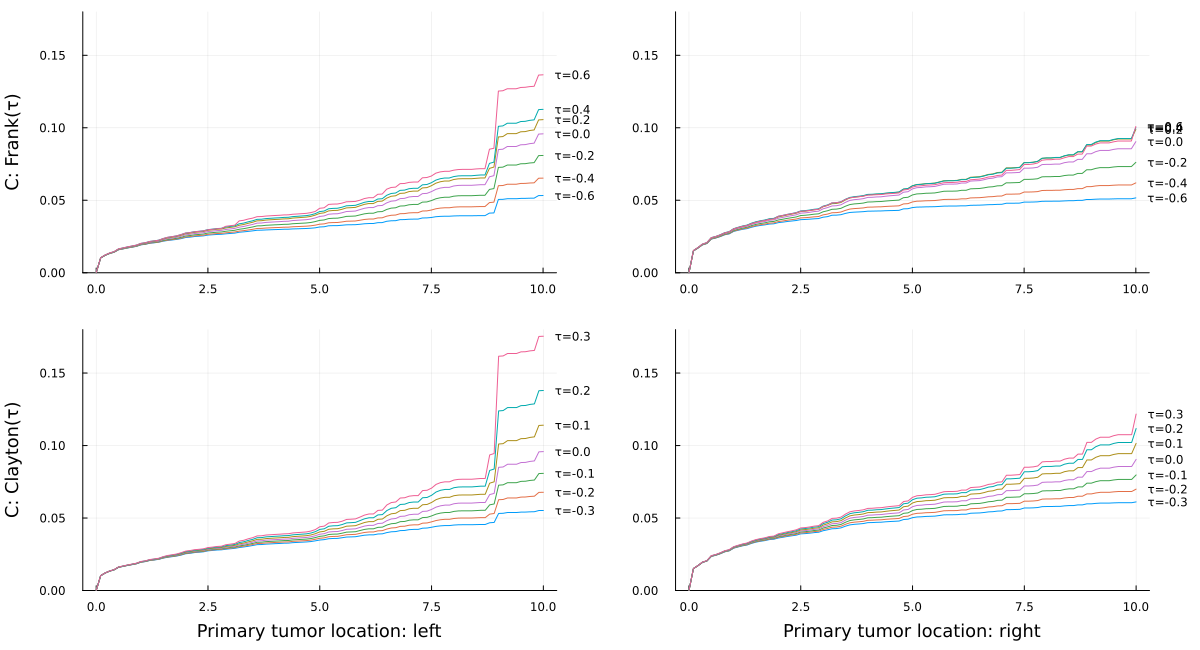}
        \caption{\label{fig:example/short_paper_graph_var} Estimated standard errors $t \to \hat{\sigma}_E(t)$ on the left and right tumor locations, for various hypothesized copulas $\mathcal C$.
        For both the Frank and Clayton copula, $\tau=0$ represents the Pohar Perme variance.}
    \end{figure}
    \end{textAtEnd}

Recall first that, on real data, the true dependence structure is unknown. 
We see on \cref{fig:example/short_paper_graph} large differences on the excess survival function with respect to the dependence structure as expected from the previous simulation.
A similar graph with the same conclusion for the variance estimate is given in \cref{apx:extra_material}.
In each graph, $\tau = 0 \iff \mathcal C = \Pi$, and this particular curve is the Pohar Perme estimator.
These results show that the selection of the dependence structure, alas a non-testable assumption, has a large impact on the practical results and conclusion a practitioner might extract from a given relative survival dataset.

\begin{figure}[H]
    \centering
    \includegraphics[width=\textwidth]{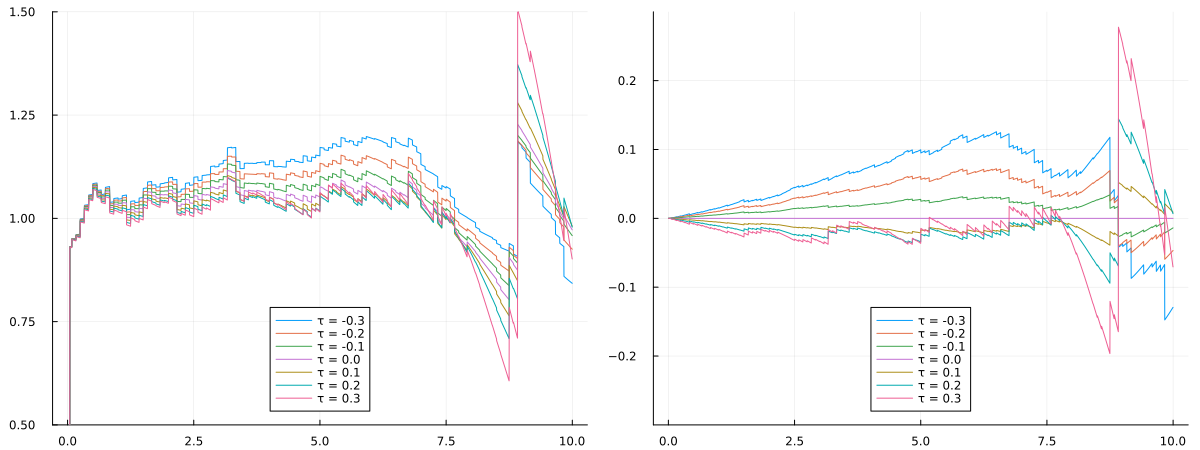}
    \caption{\label{fig:example/bootstrap_variance} The left plot represents the ratios $t \to \frac{\hat{\sigma}_E(t)}{\hat{\sigma_E}^{(\mathrm{boot})}(t)}$ for each $\texttt{Clayton}(\tau)$ copula on the left primary tumor location.
    The right plot represents the same curves, but recentered on the $\tau = 0$ one, since the variance estimator in this case (Pohar Perme) is unbiased.}
\end{figure}

To showcase on this example the potential bias in the variance estimate, we construct $N=200$ bootstrap resamples $\hat{\Lambda}_E^{(i)}, i \in 1,...,N$ of our estimator, to be able to compare the standard error estimate $\hat{\sigma}_E(t)$ from \cref{eq:gen_var} and a bootstrap version computed as 
$$\hat{\sigma_E}^{(\mathrm{boot})}(t) = \sqrt{\frac{1}{N}\sum_{i=1}^N \left(\hat{\Lambda}_E^{(i)}(t) - \frac{1}{N}\sum_{j=1}^N \hat{\Lambda}_E^{(j)}(t)\right)^2}.$$ 
The \cref{fig:example/bootstrap_variance} shows out a comparison between these two standard errors, restricted (for readability) to the case of Clayton copulas and tumors primarily located on the left side.

We see on the left plot of \cref{fig:example/bootstrap_variance} that curves are wiggling around $1$, denoting an overall low bias, but the different curves are significantly different from each others as the right plot of \cref{fig:example/bootstrap_variance} shows.
The strong effect around the $8$ years mark is due to the large estimated variance of the estimator at this time point (see \cref{fig:example/short_paper_graph_var} in \cref{apx:extra_material}).

Justification for the recentering on the Pohar Perme standard error estimates is that, under $\HPi$ (here $\tau = 0$) $\hat{\sigma}_E(t) = \tld{\sigma}_E(t)$ does not require a plug-in estimation and is thus unbiased.
We then see that, for concordant dependence structures ($\tau > 0$), $\hat{\sigma}_E(t)$'s estimated bias is negative, while it is positive for discordant dependence structures ($\tau < 0$).
Thus, we underestimate (resp overestimate) the variance when the dependence structure is concordant (resp discordant).
Hence, an unbiased bootstrap estimate of $\hat{\Lambda}_E$'s standard error would display a behavior a bit less spread out that what we already observe in \cref{fig:example/short_paper_graph_var} in \cref{apx:extra_material}.

Turning ourselves to the log-rank type tests on this dataset, the hypothesis $H_0$ is simply that the primary tumor locations (left or right) does not affect the excess survival function.
We thus run our test procedures under several $\HC$, at several ending time points $T$ (in years), to obtain p-values displayed in \cref{tab:example/tests_pvalues_clayton,tab:example/tests_pvalues_frank}.

\begin{table}[H]
    \centering\footnotesize
    \caption{\label{tab:example/tests_pvalues_clayton} Obtained p-value for the generalized log-rank-type test for $\mathcal C = \texttt{Clayton}(\tau)$, at various horizons $T$ (in years).}
    \begin{tabularx}{\linewidth}{rYYYYY}
  \toprule
  $\tau$ & $T = 2$ & $T = 3$ & $T = 5$ & $T = 8$ & $T = 10$ \\\midrule
  $-0.3$ & $0.00340$ & $0.01016$ & $0.03825$ & $0.43212$ & $0.82124$ \\
  $-0.2$ & $0.00255$ & $0.00653$ & $0.02071$ & $0.31202$ & $0.73293$ \\
  $-0.1$ & $0.00187$ & $0.00406$ & $0.01065$ & $0.21326$ & $0.61918$ \\
  $0.0$ & $0.00133$ & $0.00244$ & $0.00522$ & $0.13575$ & $0.50419$ \\
  $0.1$ & $0.00093$ & $0.00142$ & $0.00249$ & $0.08025$ & $0.40474$ \\
  $0.2$ & $0.00064$ & $0.00081$ & $0.00120$ & $0.04515$ & $0.34816$ \\
  $0.3$ & $0.00043$ & $0.00047$ & $0.00063$ & $0.02591$ & $0.41805$ \\\bottomrule
\end{tabularx}

\end{table} 

\begin{table}[H]
    \centering\footnotesize
    \caption{\label{tab:example/tests_pvalues_frank} Obtained p-value for the generalized log-rank-type test for $\mathcal C = \texttt{Frank}(\tau)$, at various horizons $T$ (in years).}
    \begin{tabularx}{\linewidth}{rYYYYY}
  \toprule
  $\tau$ & $T = 2$ & $T = 3$ & $T = 5$ & $T = 8$ & $T = 10$ \\\midrule
  $-0.6$ & $0.01208$ & $0.05266$ & $0.20128$ & $0.90222$ & $0.66067$ \\
  $-0.5$ & $0.00979$ & $0.03689$ & $0.13102$ & $0.77497$ & $0.75530$ \\
  $-0.4$ & $0.00737$ & $0.02417$ & $0.07991$ & $0.64116$ & $0.85883$ \\
  $-0.3$ & $0.00517$ & $0.01476$ & $0.04493$ & $0.49968$ & $0.98195$ \\
  $-0.2$ & $0.00342$ & $0.00845$ & $0.02329$ & $0.35883$ & $0.86804$ \\
  $-0.1$ & $0.00216$ & $0.00461$ & $0.01127$ & $0.23305$ & $0.69194$ \\
  $0.0$ & $0.00133$ & $0.00244$ & $0.00522$ & $0.13575$ & $0.50419$ \\
  $0.1$ & $0.00082$ & $0.00129$ & $0.00240$ & $0.07163$ & $0.33148$ \\
  $0.2$ & $0.00051$ & $0.00070$ & $0.00114$ & $0.03537$ & $0.19859$ \\
  $0.3$ & $0.00033$ & $0.00040$ & $0.00058$ & $0.01724$ & $0.11324$ \\
  $0.4$ & $0.00023$ & $0.00025$ & $0.00034$ & $0.00889$ & $0.06671$ \\
  $0.5$ & $0.00018$ & $0.00018$ & $0.00023$ & $0.00533$ & $0.04642$ \\
  $0.6$ & $0.00014$ & $0.00015$ & $0.00021$ & $0.00435$ & $0.04985$ \\\bottomrule
\end{tabularx}

\end{table} 

Let us first note that we clearly observe a drastic variability in the p-values with respect to the dependence structure shape and of course its Kendall's $\tau$ as we expected.
That being said, we also observe high variability with respect to  $T$.
Indeed, after $7.5$ years, the two curves have erratic behaviors and their variances drastically increase (this is visible in \cref{fig:example/short_paper_graph_var} in \cref{apx:extra_material}), making the null hypothesis much harder to reject for log-rank type tests.
However, the variance estimate being itself dependent on the assumption, a strong enough dependence structure such as the $\texttt{Frank}(\tau = 0.6)$ copula can tame this behavior down and still reject $H_0$ at these higher times.

A formal decision on whether the primary tumor location does significantly affect the net survival is therefore hard to extract from these p-values.
Focusing on one dependence structure only, e.g., historically the independence, does not seem like a good solution to this problem.
Finally, our assumption that the dependence structure is the same for each primary tumor location is not necessary, and these results could be even worse: it would indeed be possible to pick dependence structure for each group, or even each individual, to target any p-value we want.
The study on the primary tumor location thus should remain inconclusive.

\section{Discussion}\label{sec:conclusion}

The relative survival methodology is widely used by organizations and practitioners dealing with cancer registry datasets, and seeking information about the underlying cancer survival taking into account the differences in other-causes (populational) mortality of the patients.
Built on historical competing risks grounds, it relies on a hypothesis of independence between the excess mortality and the mortality w.r.t other causes, a hypothesis that is already known to be wrong but untestable with the available data.
This paper thus proposes to vary this underlying independence hypothesis through the use of copulas, in an effort to understand the consequences that this independence assumption might have on the results and conclusions drawn from the methodology.

By constructing formally a generalized version of the Pohar Perme estimator, coherent with the original version under $\HPi$, we allowed ourselves to explore the impacts of the dependency assumptions through simulations, and to perform real data analysis on a particular example.
The results are striking: the choice of the underlying dependence structure has a large impact on the resulting estimators, from the excess survival function to its variance.
Most importantly, with a wrong hypothesis for the dependence structure, the log-rank-type test becomes unusable as soon as the groups have different populational mortalities, which corresponds to most practical uses. 

The next step in the right direction would be, cancer by cancer, cohort by cohort, dataset by dataset, covariate by covariate, to try to extract this dependence structure from other datasets in which the so-called missing indicatrix is observable.
This is for example the route taken by \cite{czadoDependentCensoringBased2021} on SEER data relative to pancreas cancer, where they show a strong dependence structure between death from this cancer and censorship by other causes.
Of course, they assume the reported death cause to be true.
Without this assumption, we do not know the feasibility of such an approach.

The most pregnant issue with the proposed methodology -- and thus a relevant direction of research forward -- is the lack of analysis of the plug-in estimators, constructed from the main differential equation.
Indeed, our approach is based on the estimation of conditional probabilities $a_i,b_i,c_i$ from this differential equation, and the formal study of these estimators remains to be done.
If our experiments suggest the bias of our plugins to be negligeable, the study of the convergence of $\hat{\Lambda}_E$ to $\tld{\Lambda}_E$ (and then to $\Lambda_E$) remains open.

\appendix
\section{Proofs and secondary statements}\label{apx:proofs}
This section contains the proofs of the exposed results, preceded by a few technical lemmas upon which they stand.
\printProofs

\section{Probability distributions and copulas}\label{apx:distributions}
This section contains details on the probability distributions and copulas used in the main text.
Let us first define the following three univariate probability distributions:
\begin{itemize}
    \item $X \sim \texttt{Exponential}(\mu \in \mathbb R_+)$ if and only if $\P{X > t} = e^{-\frac{t}{\mu}}$. $\mu$ is the mean of the distribution.
    \item $X \sim \texttt{Normal}(\mu \in \mathbb R,\sigma \in \mathbb R_+)$ if and only if $\P{X \le t} = \int_{-\infty}^t \frac{1}{\sqrt{2\pi}} e^{-\frac{(x-\mu)^2}{2\sigma^2}} \partial x$.
    \item $X \sim \texttt{Chi2}(n \in \mathbb N)$ if and only if $X = \sum_{i=1}^n Z_i^2$ for $Z_i$ i.i.d $\texttt{Normal}(0,1)$.
\end{itemize}

We also extensively used the concept of copula in the paper, and some classical families of copulas which we now describe.
A copula $\mathcal C$ is essentially the distribution function of a random vector with marginals that are all uniform on $[0,1]$.
These functions are used to model dependence structure of random vectors.
We refer to standard books~\cite{cherubini2004,nelsen2006,joe2014} or more recently~\cite{durante2015a,grosser2021} on the subject.
Copulas are used to model separately the marginals and the dependence structure of random vectors, thanks to Sklar's Theorem~\ref{thm:sklar}:
\begin{theorem}[Sklar's Theorem\cite{sklar1959}]\label{thm:sklar} For every $d$-variate random vector $\bm X$ with joint distribution function $F$ and marginal distributions functions $F_i$'s, there exists a copula $C$ such that $$\forall \bm x\in \mathbb R^d, F(\bm x) = C(F_{1}(x_{1}),...,F_{d}(x_{d})).$$
  The copula $C$ is uniquely determined on $\mathrm{Ran}(F_{1}) \times ... \times \mathrm{Ran}(F_{d})$, where $\mathrm{Ran}(F_i)$ denotes the range of the function $F_i$.
  In particular, if all marginals are absolutely continuous, $C$ is unique.
\end{theorem}

There exists classical and useful parametric families of copulas.
One parametric class of interest is the class of Archimedean copulas, which can be defined in the bivariate case as follows.
A function $\phi :\mathbb R_+ \to [0,1]$ such that $\phi(0) = 1, \lim_{x\to+\infty}\phi(x) = 0$ is said to be a $2$-Archimedean generator if and only if it is $2$-monotone, i.e., non-increasing and convex.
If $\phi$ is a $2$-Archimedean generator then the function $$\mathcal C_\phi(u,v) = \phi\left(\phi^{-1}(u_1)+\phi^{-1}(u_2)\right)$$ is a bivariate copula, the so-called Archimedean copula with generator $\phi$.

Some classical generator are:
\begin{itemize}
\item $\phi(t) =e^{-t}$ which generates the independence copula $\Pi$.
\item $\phi_{\theta}(t) = \left(1+t\theta\right)^{-\theta^{-1}}$ which generates the $\texttt{Clayton}(\tau)$ copula for $\theta = \frac{2\tau}{1-\tau}$.
\item $\phi_{\theta}(t) = \exp\{-t^{\theta^{-1}}\}$ which generates the $\texttt{Gumbel}(\tau)$ copula for $\theta = \frac{1}{1-\tau}$.
\item $\phi_{\theta}(t) = -\theta^{-1}\ln\left(1+e^{-t-\theta}-e^{-t}\right)$ which generates the $\texttt{Frank}(\tau)$ copula for $\tau = 1 + \frac{4}{\theta}\int_{0}^1 \frac{t}{e^t -1}dt - \frac{4}{\theta}$ which is bijective but needs numerical integration and inversion routines.
\end{itemize}

\begin{figure}[H]
  \centering
  \includegraphics[width=\textwidth]{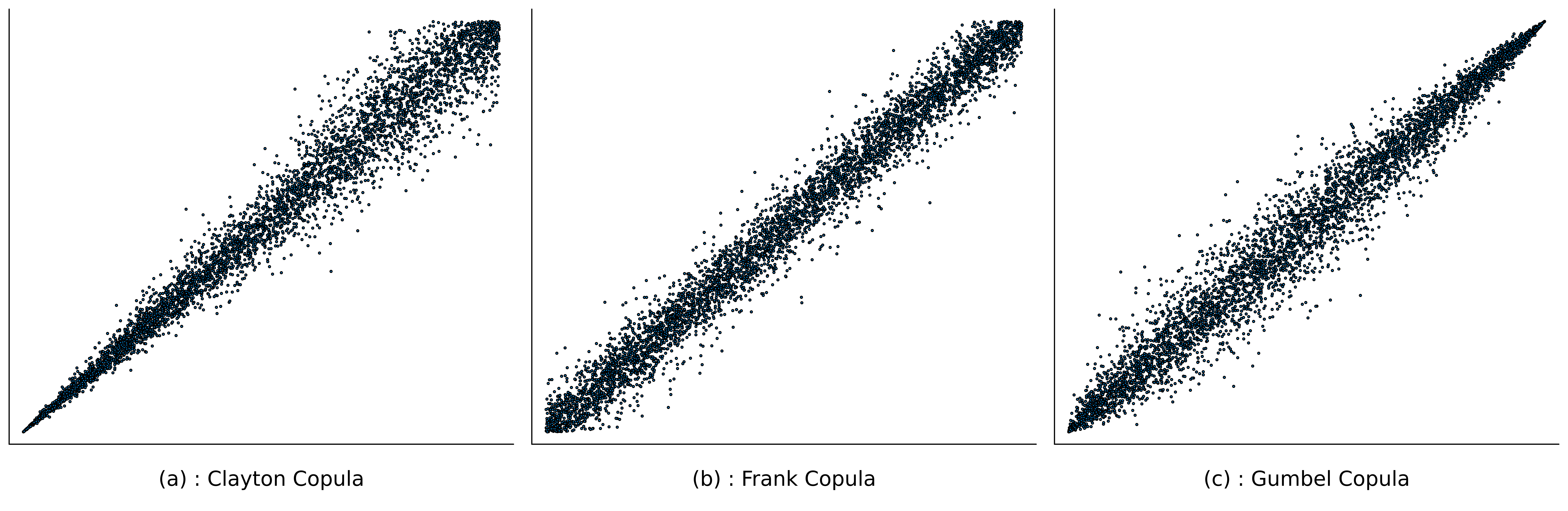}
  \caption{Simulation of $5000$ points from different Archimedean copulas, parametrized so that their Kendall's $\tau$ are all equal to $0.9$.
  Remark that the Clayton and Gumbel cases have strong dependencies in one tail only, while the Frank has both.}
  \label{fig:cops/example_archimedeans.png}
\end{figure}

We parametrized "naturally" the generators, but not the copulas.
We used a bijective mapping between generators parameters $\theta$, which have no common interpretations, and Kendall's $\tau$, which quantify the overall strength of the dependency.
Kendall's $\tau$'s take values between $-1$ and $1$ and is equal to $0$ only for the independence copula $\Pi$.
They can be computed for Archimedean copulas as:
$$\tau =1 + 4 \int_0^1 \frac{\phi^{(-1)}(u)}{\left(\phi^{(-1)}\right)'(u)} \partial u,$$ which reduces on each case to the upper conversions formulas.
Other generators are available in the proposed numerical implementation thanks to~\cite{LavernyJimenez2024}.
We represent simulated data from some example copulas in Figure~\ref{fig:cops/example_archimedeans.png}.

\section{Additional material}\label{apx:extra_material}
This section contains a few extra graphs and tables from the simulations and examples.
\printProofs[ExtraMaterial]

\printbibliography
\end{document}